
\input amstex
\documentstyle{amsppt}
\NoBlackBoxes
\magnification=1200
\nopagenumbers
\define \Fs{Fourier series}
\define \T{\bold T^n}
\define \Zn{\bold Z^n}
\define \Ff{\widehat f}
\define \ek{e^{ikx}}
\define \R{\bold R^n}

\define \Cp{{n-1 \over 2}}
\define \Ft{Fourier transform}
\define \ep{\varepsilon}
\define \Lc{Lebesgue constants}
\define \fa{\lambda}
\define \faf{\widehat \lambda}
\define \a{\alpha}
\define \ds{\partial S}
\define \ffi{\varphi}
\define \g{\gamma}

\topmatter
\author E. R. Liflyand  \endauthor
\title Estimates of Lebesgue constants \\ via Fourier transforms.
Many dimensions.     \endtitle
\affil Bar-Ilan University, Israel \endaffil
\address Department of Mathematics and Computer Science,
Bar-Ilan University, 52900 Ramat-Gan, Israel \endaddress
\email liflyand\@bimacs.cs.biu.ac.il  \endemail
\abstract This is an attempt of a comprehensive survey of the results
in which estimates of the norms of linear means of multiple Fourier
series, the Lebesgue constants, are obtained by means of estimating
the Fourier transform of a function generating such a method. Only few
proofs are given in order to illustrate a general idea of techniques
applied. Among the results are well known elsewhere as well
as less known or published in an unacceptable journals and several new
unpublished results. \endabstract
\rightheadtext{Lebesgue constants and Fourier transform}
\endtopmatter
\document
\vskip 4cm
\centerline{\bf Preprint\qquad BIMACS--9501}
\vskip 5cm
\centerline{\bf Bar-Ilan University,\qquad 1995}
\vfill \eject
\head Introduction. \endhead
{\bf 0.1.} The purpose of this work is to give a survey of results,
and partly of proofs, dealing with estimates of $L^1$-norms of linear
means of multi-dimensional Fourier series, via estimates of Fourier
transforms of functions generating these means. Such norms are called,
in a general sense, the Lebesgue constants. One might think that
Lebesgue constants in the multiple case are, in  essence, always
estimated by means of Fourier transforms. Though the answer would be
``no, not always'', such estimates are indeed the central part of
the theory, while other methods are at present only in the periphery
of the area. As for the area in general, one is reminded the words
which were said on a similar field in the marvellous book of K. M.
Davis and Y.-C. Chang [DC]: `` a topic now in disrepute due to its
difficulty''. Nevertheless, several enthusiasts continue to spent their
time seeking to find a pearl in this sea. From time to time some
intersting results appear giving hope for new growth of popularity
of this part of Fourier Analysis.

My teacher Professor R. M. Trigub seeing me returning again and again
to these problems once proposed that I should try to write a survey on Lebesgue
constants. I am very grateful to him for ``infecting'' me with this
tempting idea, and for repeatedly encouraging me to continue with the
work. Surely, there are several related surveys already in the
literature (see [Zh, AIN, AAP, Go, Dy3]), and several important
books on different aspects of Fourier Analysis; those of Stein
(see [SW, S3]) should be mentioned first. But many interesting
sides of our somewhat special topic were not touched at all in these
comprehensive works. Nevertheless, I was unsure if there was need for an
additional survey, until once Professor E. M. Stein asked me
whether it is possible to outline results on Lebesgue constants in which
the Fourier transform methods are involved. This question gave me a clear
starting point for which I would like to express my gratitude
to Professor E. M. Stein.

The outline of this survey is as follows. After preliminary notations
in \S 1 we give a proof of one V. A. Yudin's result. This proof
illustrates a general approach to such estimates. We also give several
results for the Bochner-Riesz means which were prototypes to Yudin's
result and to some other general results that are also given in this
survey. In \S 2 we give very general results mostly due to E. Belinskii.
In \S 3 we investigate various generalizations of the Bochner-Riesz
means. The next \S 4 is also devoted to generalizations of the
Bochner-Riesz means. What is preserved is the spherical nature of
summation. In \S 5 a collection of ``polyhedral'' results is given.
It is shown that there are sufficiently interesting problems in this case.
In \S 6 results are considered in which partial sums or more general
linear means are taken with respect to ``hyperbolic crosses''. In \S 7
we give as an Appendix other results which are not proved by means of
the Fourier transform methods. We tried to mention as many as possible
results on the topic; some of them are very recent.

My friends and colleagues Professors E. Belinskii, A. Podkorytov and
M. Skopina were the readers of earlier variants of this work.
Many improvements are due to their precise remarks and helpful discussions.
I would like to acknowledge their efforts and sincere interest.

Finally I wish to mention that only one person deserves ``acknowledgements''
for possible mistakes, misprints and all poor passages.
``It's Me, O Lord'', as R. Kent entitled his book.
\vfill\eject
{\bf 0.2.} Let $f$ be an integrable function on $\T$, where ${\bold T}=(-\pi,
\pi]$, $2\pi$-periodic in each variable. Consider a Fourier series of this
function        $$\sum \limits_k \Ff (k) \ek  \tag0.1 $$
where $x=(x_1,...,x_n)$ is a point in the real $n$-dimensional Euclidean
space $\R$, $k=(k_1,...,k_n) \in \Zn$, the lattice of points in $\R$
with integer coordinates, $kx=k_1 x_1 +...+ k_n x_n$ is the scalar
product, and
$$\Ff (k) = (2\pi)^{-n} \int \limits_{\T} f(x) e^{-ikx} dx$$
is the $k$-th Fourier coefficient of function $f$.

Kolmogorov's famous result \cite{K} says that the series (0.1) may be
divergent at every point of $\T$. Thus, it is quite natural to consider
a sequence of linear operators
$$L_N^\fa :\qquad f(x) \mapsto L_N^\fa (f;x) = \sum \limits_k \fa
\biggl({k \over N}\biggr) \Ff (k) \ek , \tag0.2 $$
where $\fa$ is a bounded measurable function (of course, it should be
defined at the points of type ${k \over N}$),
and to study its properties in order to derive some information about
the function $f$, or a certain space of such functions. The basic
information one receives from the behavior of the norms of these operators.
When the operators map $C(\T)$ into $C(\T)$, or $L^1 (\T)$ into $L^1 (\T)$,
and $\fa \equiv 1$ on some set and vanishes outside
(or, in the other words, $\fa$ is an indicator function
of this set) the norms are traditionally called the \Lc\. This
term is frequently saved for the general situation.

{\bf 0.3.} It is well-known \cite{SW,Ch.VII,Th.3.4} that the operator
$L_N^{\fa}$ is bounded if and only if the series
$$\sum \limits_k  \fa \biggl({k \over N}\biggr) \ek \tag0.3$$
is a Fourier series of some measure $\mu_N$, and $\Vert L_N^{\fa} \Vert
= \Vert \mu_N \Vert $. And if this series is the Fourier series of an
integrable function, the following relation takes place:
$$\Vert L_N^\fa \Vert = (2\pi)^{-n} \int \limits_{\T} \biggl| \sum \limits_k
\fa \biggl({k \over N}\biggr) \ek \biggr| dx . \tag0.4 $$
This occurs, for instance, when $\fa$ is boundedly supported.

If to take formally an integral instead of the sum on the right-hand side
of (0.4), and to fulfil some simple computations (not less formal for a
moment), considering $k$ as a continuous parameter, one can expect
something like $(2\pi)^{-n} \int_{N\T} |\faf (x)| dx$ on place of the
right-hand side of (0.4). Here
$$\faf (x) = \int \limits_{\R} \fa (u) e^{-iux} \, du $$
is the \Ft\ of $\fa$.
This process which looks so natural and so attractive (and so short!)
can be very subtle and sometimes very cumbersome in reality. We will see
below that sometimes it does not valid in some sense!

{\bf 0.4.} Very often such transformations turn on very special problems of
Number Theory, or Differential Geometry (see, e.g., the work \cite{S2} as
a representative example in a series of Stein's works devoted to the
role of the Gaussian curvature in Fourier Analysis. Besides this, in
different parts of his recent book [S3] the place of this notion among
other important things may be seen).

It is written in \cite{DC} that "the best
trick will be to transfer problems" (of convergence of Fourier series)
"from Fourier series to Fourier integrals. This is good because it is
easier to compute an integral explicitly than to sum a series in a
closed form. On the other hand, this is bad because the integrals defining
Fourier transforms do not converge absolutely".

We try to consider some problems in which "this is good" (or at least we think
so). Even in the cases when "this is bad" (in the mentioned sense),
we obtain certain information studying the order of "badness".

\head 1. Spherical partial sums, and some generalizations. \endhead
{\bf 1.1.} Let us begin with a very exemplary case. Consider the spherical
partial sums $S_N$ of the \Fs\ of a function $f$:
$$S_N (f;x) = \sum_{|k| \le N} \Ff (k) \ek ,$$
where $|k|^2 =k_1^2 + \cdots + k_n^2$. It is well-known that the norms
of operators           $$S_N :\qquad f(x) \mapsto S_N (f;x) $$
taking $C(\T)$ into $C(\T)$ (or $L^1 (\T)$ into $L^1 (\T)$ that is
the same) equal to
$$ \Vert S_N \Vert = (2\pi)^{-n} \int_{\T} \biggl|
\sum_{|k| \le N} \ek \biggr|\, dx .   \tag1.1$$
The following ordinal estimate holds.
\proclaim{Theorem 1.1} There exist positive constants $C_1$, and $C_2$,
where $C_1 < C_2$, depending only on $n$, such that
$$C_1 N^{\Cp} \le \Vert S_N \Vert \le C_2 N^{\Cp} .   \tag1.2$$
\endproclaim
\demo{Proof} We give an outline of the proof which illustrates a method,
rather general and having the passage to the \Ft\ described in the
Introduction as the basic construction.

The estimate from below was obtained firstly by V.A.Ilyin \cite{I} (for
expansions corresponding to the Laplace operator), and after that two-sided
estimates were obtained in \cite{Ba2} and \cite{IA}. All those proofs
were rather complicate. We will follow a very simple proof of the
upper estimate proposed by Yudin \cite{Y1} (in more
general situation; the earlier proof of the upper estimate in (1.2) due to
H.Shapiro \cite{Sh} essentially is almost the same). Then, using
one trick proposed in \cite{I}, we will adjust this proof for the
estimate from below as well.

If $I_k$ is the cube with the edge of length $1$ and the center at the point
$k$ , and $B_N = \bigcup_{k: |k| \le N} I_k$, then
$$\align  \int\limits_{B_N} e^{iux} du &= \sum_{k: |k| \le N}
\int\limits_{I_k} e^{iux} du \\ \quad\\
= \sum_{k: |k| \le N} \ek  \prod \limits_{j=1}^n {2 \sin {x_j \over 2}
\over x_j} &= \prod \limits_{j=1}^n {2 \sin {x_j \over 2} \over x_j}
\sum_{k: |k| \le N} \ek .  \endalign$$
Thus, we obtain
$$\align  \Vert S_N \Vert &\le (2\pi)^{-n} ({\pi \over 2})^n \int\limits_{\T}
\biggl|\quad \int\limits_{B_N} e^{iux} du \biggr|\, dx \\ \quad\\
&\le 4^{-n}\int\limits_{\T} \biggl|\quad \int\limits_{|u| \le N} e^{iux} du
\biggr|\, dx  + 4^{-n} \int\limits_{\T} \biggl|\quad
\int\limits_{D_N} e^{iux} du \biggr|\, dx ,  \tag1.3 \endalign$$
\medskip\flushpar
where $D_N$ is a symmetric difference of the sets $B_N$ and $|u| \le N$.
Taking into account that $\operatorname{mes} D_N \le CN^{n-1}$ and
applying the Cauchy-Schwarz inequality to the last summand
on the right-hand side of (1.3), we obtain by virtue of Parseval's
equality that
$$\align \Vert S_N \Vert & \le 4^{-n} \int\limits_{\T} \biggl|\ \int\limits_
{|u|\le N} e^{iux} du \biggr|\, dx + O(N^{\Cp}) \\ \quad\\
& = 4^{-n} \int\limits_{N\T} \biggl|\ \int\limits_{\R} \chi_1 (u) e^{iux} du
\biggr|\, dx + O(N^{\Cp}) \\ \quad\\
& = 4^{-n} \int\limits_{N\T} |\hat \chi_1 (x) | \, dx + O(N^{\Cp}) ,
\tag1.4 \endalign$$
where $\chi_1$ is the indicator function of the unit ball $|u| \le 1$.
Now standard computations of $\hat \chi_1$ via Bessel functions and
consequent integration (see, e.g., \cite{SW}, Ch.7) complete the upper
estimate in (1.2).

In order to obtain the lower estimate, let us introduce a small
parameter $\ep$, $0 < \ep < 1$, as it was done in \cite{I}. We obtain,
instead of (1.1),
$$ \Vert S_N \Vert \ge (2\pi)^{-n} \int \limits_{\ep \T}
\biggl| \sum_{|k| \le N} \ek \biggr| \, dx . $$
Now we can repeat the operations similar to (1.3), (1.4), but with
estimates from below and signs "$-$" instead of "+" on appropriate
places. This will give us
$$\Vert S_N \Vert \ge (2\pi)^{-n} \int\limits_{\ep N\T} |\hat \chi_1 (x)| \, dx
- C_3 \ep^{{n \over 2}} N^{\Cp} .  \tag1.5$$
The afore-mentioned computations via Bessel functions will yield here
$$\Vert S_N \Vert \ge \left( C_4 \ep^{\Cp} - C_3 \ep^{{n \over 2}}
\right) N^{\Cp} .$$
It remains only to choose $\ep$ such that $C_1 =
C_4 \ep^{\Cp} - C_3 \ep^{{n \over 2}} > 0$.
The proof is complete. \qquad\qed   \enddemo

In the less known paper by A. Podkorytov [P0] similar technique was
elaborated independently. This allowed to obtain the following
interesting result. To indicate that partial sums correspond to
certain set $B$ we will denote them by
$$S_B(f;x)=\sum\limits_{k\in B}\hat f(k) e^{ikx} .\tag1.6$$
\proclaim{Theorem 1.2} Assume that the set $B\subset\Bbb R^n$ satisfies
the following conditions:
\medskip
1) $B\subset[-N_1,N_1; -N_2,N_2;...;-N_n,N_n],$ where $N_1\ge
N_2\ge...\ge N_n\ge 0.$
\medskip
2) For all $j=1,...,n$ and all $x_1,...,x_{j-1},x_{j+1},...,x_n$ the set
$$B_j(x_1,...,x_{j-1},x_{j+1},...,x_n)=\{x_j: (x_1,...,x_j,...,x_n)\in B\}$$
is either empty or is an interval.

Then $$||S_B||=O\biggl(\sqrt{N_2\cdots N_n}\biggl(1+\ln{N_1\over
N_2}\biggr)\biggr).$$          \endproclaim

{\bf 1.2} This way, or certain of its modifications, is used in many results
dealing with the estimates of \Lc\ of partial sums as well as the
linear means of \Fs\. We have already mentioned that such was  Yudin's
estimate from above \cite{Y1}, for more general sets $B$ generating the
corresponding partial sums, namely those which are balanced (with each
point $x$ the whole set $\delta x$, $|\delta| \le 1$, belongs to $B$),
and having the finite upper Minkowski measure:
$$\limsup_{\ep \to 0} {1 \over \ep} \operatorname{mes}
\{ x: \rho (x, \partial B) < \ep \} < \infty , $$
where $\rho (x,y)$ is the distance between two points $x, y \in \R$,
and $\rho (x, \partial B) = \inf \limits_{y \in \partial B} \rho (x,y)$.

The same method was applied for estimate from below in \cite{L2},
where conditions are less restrictive (above all, they are local)
than in the earlier paper \cite{CaS} and the later papers [Br1, Br2].

\proclaim{Theorem 1.3} {\rm ([L2, L3])} Let the boundary of the
region $B$ contain a simple (non-intersecting) piece of a surface
of smoothness $[{n+2 \over 2}]$ in which there is at least one point with
non-vanishing principal curvatures. Then there exists a positive
constant $C$ depending only on $B$ such that
$$\int\limits_{\T} \left| \sum_{k \in N\!B \cap \Zn} \ek \right| \, dx
\ge C N^{\Cp} $$     for large $N$.         \endproclaim \medskip
We will see below a generalization of this theorem (see Theorem 3.2).
And now let us give one related two-dimensional result.
\proclaim{Theorem 1.4} {\rm ([Gu])} Assume that a convex set $B$ is
included into $\Bbb T^2.$ Then for sufficiently large $N$ the inequality
$$||S_{NB}||\ge CN^{{1\over 2}}\biggl(\int\limits_{-\pi}^\pi\sqrt{\rho
(\varphi)}\,d\varphi\biggr)^2$$ holds, where $\rho(\varphi)$ is the
curvature radius of $\partial B$ at the point where $\max\{x_1\cos\varphi
+x_2\sin\varphi: x\in B\}$ is attained.

If $\liminf{||S_{NB}||\over N^{{1\over 2}}}=0$ then the boundary of
the set $B$ is degenerate, i.e. for almost all directions $\varphi$
its curvature radius is equal to zero.  \endproclaim
What is of special interest in this result, as well as in [P7],
is the fact that no assumptions on smoothness of $\partial B$ are
involved. Of course, convexity itself gives some minimal smoothness.
Then, only two-dimensional statements are obtained in [P7] and [Gu].
This gives rise to the question: what are minimal smoothness
assumptions for such estimates in case of arbitrary dimension.
\bigskip
\head 2. Some general estimates of Lebesgue constants. \endhead
{\bf 2.1.} Certainly, we were speaking about an outline, and in more general
and more complicate situations each step can be cumbersome and
entailed with bigger technical difficulties. For example, even for the
case of Fejer means of partial sums generated by convex sets, a proof
of the boundedness of the norms of corresponding operators in \cite{P1}
is rather complicate, and all the time refers to non-trivial
estimates of \Ft s.

{\bf 2.2.} Let us give now one Belinskii's result in which this method is
realized on a very high level of generality. Belinskii was apparently the first
began a {\it systematic\/} study of connections between summability and
integrability of the Fourier transform of a function generating a
method of summability, in the multi-dimensional case.
\proclaim{Theorem 2.1} {\rm ([Be2])} Let $\fa$ be a bounded
measurable function with a compact support. Then for the norms of a
sequence of linear operators {\rm (0.2)} we have
$$\align \Vert L_N^{\fa} \Vert_{L^1 (\T) \to L^1 (\T)} & \le
(2\pi)^{-n} \int \limits_{N\T} \prod \limits_{j=1}^n {x_j \over
2N \sin {x_j \over 2N}} |\faf (x)| \, dx \\
& + \sum_{j=1}^{m-1} ({\pi \over 2})^{(j+1)n} \int \limits_{N\T} |\faf (x)|
\left|{x \over N}\right|^j \, dx \tag2.1 \\
+ {\pi^{nm+n/2} \over 2^{nm-n/2}} \int \limits_{{1 \over 2\pi}\T} & \cdots
\int \limits_{{1 \over 2\pi}\T} \left( \sum_k \left| \Delta_{{k \over N}}^m
\left( \fa ; {u_1 \over N},...,{u_m \over N} \right) \right|^2
\right)^{{1 \over 2}} du_1 ... du_m ,  \endalign$$

\quad

$$C_p \Vert L_N^{\fa} \Vert_{L^p (\T) \to L^p (\T)}  \ge\left\{(2\pi)^{-n}
\int\limits_{\ep N\T} \left| \prod \limits_{j=1}^n {x_j \over
2N \sin {x_j \over 2N}} \faf (x) \right|^p \, dx \right\}^{{1 \over p}} $$
$$- \sum_{j=1}^{m-1} ({\pi \over 2})^{(j+1)n} \left\{ \int \limits_{\ep N\T}
|\faf (x)|^p \left|{x \over N}\right|^{jp} \, dx
\right\}^{{1 \over p}} \tag2.2 $$
$$- \dfrac{\pi^{nm+n{2-p \over 2p}}}{2^{nm-n{2-p \over 2p}}}
\frac{\ep^{n{2-p \over 2p}}}{N^{n{p-1 \over p}}}
\int \limits_{{1 \over 2\pi}\T} \cdots
\int \limits_{{1 \over 2\pi}\T} \left( \sum_k \left| \Delta_{{k \over N}}^m
\left( \fa ; {u_1 \over N},...,{u_m \over N} \right) \right|^2
\right)^{{1 \over 2}} du_1 ... du_m .$$
\endproclaim

\quad

Here $\ep$, $0<\ep \le 1$, is an arbitrary real number, $m$ is integer, and
$1 \le p \le 2.$ The $m$-th difference $\Delta_z^m (\fa; h_1 ,..., h_m)$
is defined recursively by the formulas
$$\Delta_z^1 (\fa ; h_1) = \fa (z+h_1) - \fa (z) ;$$
$$\Delta_z^m (\fa; h_1 ,..., h_m) = \Delta_{z+h_m}^{m-1}
(\fa; h_1 ,..., h_{m-1}) - \Delta_z^{m-1} (\fa; h_1 ,..., h_{m-1}) ,$$
with $h_j , z \in \R$. When $p > 2$, in view of duality (see, e.g.,
\cite{SW,Ch.I,Th.3.20}), the estimate (2.2) still
valids with $p' = {p \over p-1}$ instead of $p$.
We give the proof of the main case: the estimate from below for $p=1.$
\demo{Proof} We follow the argument from \cite{Be2}. We have by definition
$$\Vert L_N^{\fa} \Vert_{L^p (\T) \to L^p (\T)}=\sup\limits_{||f||_
{L^p(\T)}\le 1}\left\{(2\pi)^{-2n}\int\limits_{\T}\left|f(x-u)
\sum\limits_k \fa \biggl({k\over N}\biggr) e^{ik x}
\right|^p \,dx\right\}^{{1\over p}}.$$
For $p=1$ we have
$$||L_N^{\fa}||=(2\pi)^{-n} \int\limits_{\T}
\biggl|\sum\limits_k \fa \biggl({k\over N}\biggr) e^{ik x}
\biggr|\,dx.\tag0.4$$
For $p>1$, assume without loss of generality that $\fa=0$ when
$|x_j|>1,$ $j=1,2,...,n,$ and set
$$f(x_1,...,x_n)={N^{n{1-p\over p}}\over C_p}\prod\limits_{j=1}^n
D_N (x_j)$$ where $D_N(x_j)$ is the Dirichlet kernel. A constant
$C_p$ is chosen to provide $||f||_{L^p(\T)}\le 1.$ We obtain
$$C_p \Vert L_N^{\fa} \Vert_{L^p (\T) \to L^p (\T)}\ge N^{n{1-p\over p}}
\left\{(2\pi)^{-n} \int\limits_{\T}
\biggl|\sum\limits_k \fa \biggl({k\over N}\biggr) e^{ik x}
\biggr|^p \,dx\right\}^{{1\over p}}.$$
Using the obvious equality $$\int\limits_{k_j -{1\over 2}}^{k_j
+{1\over 2}} e^{ix_j \cdot v_j}\, dv_j = e^{ik_j\cdot x_j} {2\sin
{x_j \over 2}\over x_j} , \qquad j=1,2,...,n,$$
replace the sum by the Fourier transform of $\fa$ as in \cite{Y1}
(see also \cite{Zg}, Ch.V, Th.2.29):
$$\align & C_p \Vert L_N^{\fa} \Vert_{L^p (\T) \to L^p (\T)} \\
\qquad\\  &=N^{n{1-p\over p}}\left\{(2\pi)^{-n} \int\limits_{\T}
\biggl|\prod\limits_{j=1}^n {x_j\over 2\sin{x_j\over 2}}
\sum\limits_k \int\limits_{k+{1\over 2\pi}\T}
\fa\biggl({k\over N}\biggr) e^{ixv}\,dv\biggr|^p\,dx
\right\}^{{1\over p}}. \endalign$$
Obviously that for $0<\ep <1$ \medskip
$$\align & C_p \Vert L_N^{\fa} \Vert_{L^p (\T) \to L^p (\T)} \\ \quad\\&\ge
N^{n{1-p\over p}}\left\{(2\pi)^{-n} \int\limits_{\ep\T}
\biggl|\prod\limits_{j=1}^n {x_j\over 2\sin{x_j\over 2}}
\sum\limits_k \int\limits_{k+{1\over 2\pi}\T}
\fa\biggl({k\over N}\biggr) e^{ixv}\,dv\biggr|^p\,dx
\right\}^{{1\over p}}. \endalign$$ \medskip
The simple inequality ${2\over\pi}t\le\sin t ,$ whis $0\le t\le{\pi\over 2},$
and Minkowski's inequality yield \medskip
$$\align & C_p \Vert L_N^{\fa} \Vert_{L^p (\T) \to L^p (\T)} \\ \quad\\
& \ge N^{n{1-p\over p}}\left\{(2\pi)^{-n} \int\limits_{\ep\T}
\biggl|\prod\limits_{j=1}^n {x_j\over 2\sin{x_j\over 2}}
\sum\limits_k \int\limits_{k+{1\over 2\pi}\Bbb T^n}
\fa\biggl({v\over N}\biggr) e^{ixv}\,dv\biggr|^p\,dx\right\}^{{1\over p}} \\
\quad\\ & -N^{n{1-p\over p}}\left\{(2\pi)^{-n} ({\pi\over 2})^{np}
\int\limits_{\ep\T}\biggl|\sum\limits_k \int\limits_{k+{1\over 2\pi}\T}
\biggl[\fa\biggl({k\over N}\biggr) -\fa\biggl({v\over N}\biggr)
\biggr] e^{ix\cdot v}\,dv\biggr|^p\,dx\right\}^{{1\over p}}.\endalign$$
Summing and change of variables reduces the first term on the right-hand
side to the form claimed. Let us estimate the second term. Change of
variables $v_j \to v_j +k_j ,$ for $j=1,...,n,$ and generalized Minkowski's
inequality yield
$$\align &N^{n{1-p\over p}}({\pi\over 2})^n\left\{(2\pi)^{-n}\int\limits_
{\ep\T}\biggl|\sum\limits_k\int\limits_{k+{1\over 2\pi}\T}\biggl[\fa
\biggl({k\over N}\biggr)-\fa \biggl({v\over N}\biggr)\biggr]
e^{ixv}\,dv\biggr|^p\,dx\right\}^{{1\over p}} \\ \quad\\
& \le N^{n{1-p\over p}}(2\pi)^{-{n\over p}} ({\pi\over 2})^n\int
\limits_{{1\over 2\pi}\T}\left\{\int\limits_{\ep\T}\biggl|\sum\limits_k
\biggl[\fa\biggl({k\over N}\biggr)-\fa \biggl({k+v\over N}\biggr)\biggr]
e^{ikx}\biggr|^p\,dx\right\}^{{1\over p}}\,dv \\ \quad\\
&=N^{n{1-p\over p}}(2\pi)^{-{n\over p}}({\pi\over 2})^n
\int\limits_{{1\over 2\pi}\T}\left\{\int\limits_{\ep\T}\biggl|
\sum\limits_k \Delta_{{k\over N}}^1(\fa ;{u_1\over N}) e^{ikx}
\biggr|^p\,dx\right\}^{{1\over p}}\,du_1 .\endalign$$ \medskip
Note that the inner integral on the right-hand side is of the same form
as in the beginning of the proof, so the same argument is
applicable to it. In what follows we need only estimates
from above. Taking into account that \medskip
$$\align & \int\limits_{\R}  \biggl[\fa\biggl({u_2\over N}\biggr)
-\fa\biggl(u_2 +{u_1\over N}\biggr)\biggr] e^{-iu_2 x}\,du_2 \\ \quad\\
=&\int\limits_{\R}\fa (u_2)[e^{-iu_2 x}-e^{-i(u_2 -{u_1\over N}) x}]\,du_2 \\
\quad\\ =[1-e^{i{u_1\over N} x}]&\int\limits_{\R}\fa(u_2)e^{-iu_2 x}
\,du_2 =[1-e^{i{u_1\over N} x}]\widehat\fa (x) ,\endalign$$
we obtain
$$\align & N^{n{1-p\over p}}(2\pi)^{-{n\over p}}({\pi\over2})^n\int\limits_
{{1\over 2\pi}\T}\left\{\int\limits_{\ep\T}\biggl|\sum\limits_k\Delta_{{k\over
 N}}^1(\fa ;{u_1\over N})e^{ikx}\biggr|^p\,dx\right\}^{{1\over p}}\,du_1\\
& \le N^{n{1-p\over p}}(2\pi)^{-{n\over p}}({\pi\over 2})^{2n}\int\limits_
{{1\over 2\pi}\T}\left\{\int\limits_{\ep N\T} \biggl|[1-e^{i{u_1\over N} x}]
\widehat\fa (x)\biggr|^p\,dx\right\}^{{1\over p}}\,du_1 \\
&+N^{n{1-p\over p}}(2\pi)^{-{n\over p}}({\pi\over 2})^{2n}\int\limits_{{1\over
 2\pi}\T}\int\limits_{{1\over 2\pi}\T}\left\{
\int\limits_{\ep\T} \biggl|\sum\limits_k \Delta_{{k\over N}}^2
(\fa ;{u_1\over N},{u_2\over N}) e^{ikx}\biggr|^p\,dx\right\}^{{1\over p}}
\,du_1\,du_2 .  \endalign$$
Since $\sin t \le t$ for $t>0 ,$ we get using the Cauchy-Schwarz inequality
$$|1-e^{i{u_1\over N} x}| \le {|u_1|\over N} |x| .$$  Hence
$$\int\limits_{{1\over 2\pi}\T}\left\{\int\limits_{\ep\T} \biggl|
[1-e^{i{u_1\over N}x}]\widehat\fa (x)\biggr|^p\,dx\right\}^{{1\over p}}\,du_1
\le\left\{\int\limits_{\ep N\T}\biggl|{x\over N}\biggr|^p |\widehat\fa(x)|
\,dx\right\}^{{1\over p}} .$$
Repeating the same computations $m-2$ times, we obtain (2.2) with
the remainder term
$$N^{n{1-p\over p}}(2\pi)^{-{n\over p}}({\pi\over 2})^{nm}\int\limits_{{1\over
 2\pi}\T}... \int\limits_{{1\over 2\pi}\T}\left\{\int\limits_{\ep\T}
\biggl|\sum\limits_k \Delta_{{k\over N}}^m (\fa ;{u_1\over N}
,...,{u_m\over N}) e^{ikx}\biggr|\,dx\right\}^{{1\over p}}\,du_1 ... du_m .$$
In order to complete the proof apply H\"older's inequality with the power
${2\over p}$ to the inner integral, and then Parseval's equality.
The opposite inequality ($p=1$) can be obtained similarly. \qquad\qed  \enddemo

{\bf 2.3.} To show the strength of this theorem, we can say that not only
Theorem 1.1 follows from it as a technical corollary, but the following more
general result as well.

Let $\fa (x)= R_{\a} = (1-|x|^2)^{\a}_+$. Corresponding linear means
generated by this function are called the {\it Bochner-Riesz means}. An
unbreakable interest to these means was initiated by Bochner's famous
work \cite{Bc}.
\medskip
\proclaim{Theorem 2.2} {\rm ([I, Ba2, IA])}
There exist positive constants $C_1$ and $C_2$, where $C_1 < C_2$, depending
only on $n$ and $\a$, such that for $0 \le \a < \Cp$
$$C_1 N^{\Cp - \a} \le \left \Vert L_N^{R_{\a}} \right \Vert_{L^1
(\T) \to L^1 (\T)} \le C_2 N^{\Cp - \a} .   \tag2.3$$
\endproclaim
\bigskip
Moreover the following estimate due to Babenko \cite{Ba2}, \cite{Ba3}
follows from Theorem 2.1 in the same manner:
\proclaim{Corollary 2.1} The following estimate holds:
$$\left \Vert L_N^{R_{\a}} \right \Vert_{L^p (\T) \to L^p (\T)} \ge C
\left(\frac{N^{p(\a_p - \a)}-1}{p(\a_p - \a)} \right)^{{1 \over p}} ,\tag2.4$$
where $1 \le p \le {2n \over n+1}$, and $0 \le \a < \a_p =
{n \over p} - {n+1 \over 2}$. For $\a = \a_p$ this estimate should be
understood as a limit as $\a \to \a_p$, and gives a logarithmic
order of growth. \endproclaim
\demo{Proof} The Fourier transform of $R_{\a}$ is very well known (see
e.g., [SW], Ch.4):
$$\widehat R^{\a}(x)=2^{-{n\over 2}+\a}\pi^{-{n\over 2}}\Gamma(\a+1)
|x|^{-{n\over 2}-\a} J_{{n\over 2}+\a}(|x|)$$
where $J_\nu$ is the Bessel function of the first kind and order
$\nu.$ Let us estimate in (2.1) summands with the Fourier transform.
In each one extend the domain of integration to the ball of the radius
$\sqrt{n}\pi N$ and pass to spherical coordinates. We get
$$||L_N^{R_{\a}}||\le C\int\limits_0^{\sqrt{n}\pi N}\biggl| {J_{{n\over 2}+\a}
(r)\over r^{{n\over 2}+\a}}\biggr|r^{n-1}\,dr+R,$$
where $R$ denotes the last summand in (2.1) or (2.2). We will estimate
it separately. Let us use the following asymptotic formulas for Bessel
functions (see e.g., [BE], 7.12(8),7.13.1(3)):
$$J_\nu(t)=\biggl({t\over 2}\biggr)^\nu{1\over\Gamma(\nu+1)}+
O(|t|^{\nu+2}),\qquad\text{as}\quad t\to 0. \tag2.5$$
$$J_\nu(t)=\biggl({2\over\pi t}\biggr)^{{1\over 2}}\cos(t-{\nu\pi\over 4}
-{\pi\over 4})+O(t^{-{3\over 2}}),\qquad\text{as}\quad t\to\infty.\tag2.6$$
These yield $$\int\limits_0^{\sqrt{n}\pi N}\biggl| {J_{{n\over 2}+\a}
(r)\over r^{{n\over 2}+\a}}\biggr|r^{n-1}\,dr\le C{N^{{n-1\over 2}-\a}-1\over
{n-1\over 2}-\a}+O(N^{{n-3\over 2}}).\tag2.7$$
Let us estimate now the remainder $R$.
$$\align &\sup\limits_{u_1,...,u_m\in{1\over 2\pi}\T}\biggl\{\sum\limits_k
\biggl|\Delta_{{k\over N}}^m\biggl(R_{\a};{u_1\over N},...,{u_m\over N}
\biggr)\biggr|^2\biggr\}^{{1\over 2}}\\ \quad\\
=& \sup\limits_{u_1,...,u_m\in{1\over 2\pi}\T}\biggl\{\sum\limits_{|k|<
N-m-2}+\sum\limits_{N-m-2\le |k|\le N}\quad \biggr\}^{{1\over 2}}.\endalign$$
Estimate each summand in the second sum by maximal value of the function
at these points. Since $R_{\a}$ is monotone increasing near the origin
we have $$\sum\limits_{N-m-2\le |k|\le N}\le CN^{-2\a} \sum\limits_{N-m-2\le
|k|\le N} 1\le CN^{n-1-2\a}$$ (the latter value follows from the well known
estimates of the number of points of $\Bbb Z^n$ in the $n$-dimensional
ball of radius $N.$ The mean-value theorem for the directional derivative
yields that the first sum is
$$N^{-2m}\sum\limits_{|k|< N-m-2}\biggl|{\partial^m R_{\a}\over
\partial u_1...\partial u_m} \left.\right|_{{k\over N}+\theta_1{u_1\over N}
+...+\theta_m{u_m\over N}}\biggr|^2.$$
Choose $m$ such that $m\ge\a+1.$ If $\a$ is integer than the derivative
is bounded and $$\sum\limits_{|k|<N-m-2}\le CN^{n-2-2\a}.$$
If $\a$ is fractional than we have by estimating the derivative by its
maximal value on the interval:
$$\sum\limits_{|k|<N-m-2}\le CN^{-2m}\sum\limits_{|k|<N-m-2}
\biggl(1-{(|k|+m+2)^2\over N^2}\biggr)^{2(\a-m)}.$$
In view of monotonicity of the function it is possible to integrate in
place of summing. This gives $$\sum\limits_{|k|<N-m-2}\le CN^{n-1-2\a},$$
and finally $$R\le CN^{{n-1\over 2}-\a}.\tag2.8$$
Use now (2.5), (2.6), and (2.8). The formula (2.2) yields
$$||L_N^{R_{\a}}||_{L^p(\T)\to L^p(\T)}\ge C_1\biggl[{N^{p(\a_p-\a)}-1
\over p(\a_p-\a)}\biggr]^{{1\over p}}\ep^{\a_p-\a}-C_2 N^{\a_p-\a}
\ep^{\a_p+{1\over 2}}.$$ Choose $\ep$ so that
$$C_1\ep^{\a_p-\a}-C_2\ep^{\a_p+{1\over 2}}\ge C>0,$$
and this completes the proof.  \qquad\qed\enddemo
\remark{Remark 2.1} Observe that (2.7) and (2.8) give the right-hand
side of (2.3). \endremark
The original proofs and the way of derivation them from
Theorem 2.1 are uncomparable as for hardnesss. For instance,
properties of Riemann's Zeta Function is the main tool in [Ba2].
It is more convenient for us to give Corollary 2.2 in the next
section.

{\bf 2.4.} Let us give some other results of Belinskii \cite{Be1}.
These results are based on the Poisson summation formula and
technique of estimating trigonometric sums and integrals via the \Ft\.

Consider a certain function $\fa (x)$ bounded on $\R$ and continuous
at the points of $\Zn$, and construct the formal trigonometric series
for $f \in C(\T)$
$$\sum_{k \in \Zn} \fa (k) \Ff (k) \ek .\tag2.9$$
Denote $$U_N (f;x) = (2\pi)^{-n} \int\limits_{\T} f(x+u)\sum_m \fa
\biggl({m \over N}\biggr) e^{-imu}\, du . $$
The following two propositions are very useful in many applications.
They are the corollaries to Theorem 1 in \cite{Be1}. We do not formulate
the theorem itself because just these propositions have proved to
concentrate its main possibilities.
\proclaim{Proposition 2.1} Suppose $U_N (f;x)$ is defined as above.
If $\fa \in C(\R)$ and $\faf \in L^1 (\R)$, then
$$\Vert U_N \Vert = (2\pi)^{-n} \int\limits_{\R} |\faf (u)|\,du +
\theta (2\pi)^{-n} \int\limits_{\R \setminus N\T} |\faf (u)| \, du , $$
where $-2 \le \theta \le 0$, and, in general, the constant $\theta$
depends on $N$.               \endproclaim
\bigskip
\proclaim{Proposition 2.2} If $\fa$ is boundedly supported and continuous,
then
$$\sup \limits_N \Vert U_N \Vert = (2\pi)^{-n} \int\limits_{\R} |\faf (u)|\,
du .  \tag2.10$$              \endproclaim
\bigskip
All the norms here are the uniform norms.
Note, that similar upper estimate, in the case when both $\fa$ and $\faf$
are integrable, may be found in \cite{SW, Ch.VII, Sect.2}:
$$\int \limits_{\T} \biggl| \sum \limits_k \fa \biggl({k \over N}\biggr) \ek
\biggr| dx \le \int \limits_{\R} | \faf (x) | dx .  \tag2.11$$
\bigskip
{\bf 2.5.} In many questions dealing with summability some assumptions on
$\fa$, connected with bounded variation, are rather natural. We can notice
several such works in the one-dimensional case, say, [H, Te, Be0, T2].
Let us give one result of Trigub,
which is quite general and seems to be very useful in many cases when
one needs to pass from Fourier series to Fourier integrals. We need
some well-known notions. The Vitali variation is defined as follows
(see e.g., \cite{AC}). Let $\fa$ be a complex-valued function and
$$\Delta_u \fa (x)=\biggl(\prod\limits_{j=1}^n \Delta_{u_j}\biggr)\fa (x) , $$
$$\Delta_{u_j}\fa(x)=\fa(x)-\fa(x_1,...,x_{j-1},x_j +u_j ,x_{j+1} ,...,x_n),$$
be a ``mixed'' difference with respect to parallelepiped $[x,x+u]$.
Let us take an arbitrary number of non-overlapping parallelepipeds, and form
a mixed difference with respect to each of them. Then the Vitali variation is
$$V(\fa) = \sup \sum |\Delta_u \fa (x)| $$
where the least upper bound is taken over all the sets of such
parallelepipeds. For smooth functions $\fa$
$$V(\fa) = \int \limits_{\R} \left | \frac{\partial^n \fa (x)}
{\partial x_1 ... \partial x_n} \right | \, dx .$$
The Tonelli variation is something another \cite{To}. Roughly saying,
a function is of bounded Tonelli variation if it has a bounded variation
in each variable, and these variations are integrable as functions of
the rest variables. For smooth function $\fa$ it is equal to
$$\int\limits_{\R} \sum_{j=1}^n \left|{\partial \fa (x) \over
\partial x_j } \right| \, dx .$$
Let us write $\fa \in V_0$ if its Vitali variation is bounded and
$\lim \limits_{|x| \to \infty} \fa (x) =0$. In this case the function
is of bounded variation with respect to any smaller number of variables.
\proclaim{Theorem 2.3} {\rm ([T2, T3])} The following relations hold:

a) For each $\fa \in V_0$, and for every $\ep = (\ep_1,...,\ep_n)$,
$\ep_j > 0$, $j=1,...,n$,
$$ \sup \limits_{0 < |u_j| \le {\pi \over \ep_j}}  \biggl|
\int\limits_{\R} \fa (x) e^{-iux} \, dx  -  \prod \limits_{j=1}^n \ep_j
\sum \limits_k \fa (\ep_1 k_1 ,..., \ep_n k_n ) e^{-i(\ep_1 k_1 u_1 +
... + \ep_n k_n u_n )} \biggr| $$
$$ \le C V(\fa) \sum \limits_{j=1}^n \ep_j \prod \limits_{q \not = j}
{1 \over |u_q|} .  \tag2.12 $$

b) If, moreover, $\fa$ has also a bounded Tonelli variation,
dominated by $V(\fa)$, then in {\rm (2.12)} it is possible to put ${1 \over
1+|u_q|}$ instead of ${1 \over |u_q|}$.
\medskip
c) If $\fa$ satisfies a) and b) then for $N=(N_1,...,N_n)$ and
${k \over N} = ({k_1 \over N_1},...,{k_n \over N_n})$
$$\Vert L_N^{\fa} \Vert = (2\pi)^{-n} \int \limits_{|x_j| \le \pi N_j}
|\faf (x)|\, dx + \theta V(\fa) \sum \limits_{j=1}^n \prod \limits_{q
\not = j} \ln (N_q +1), \quad |\theta|\le C . \tag2.13 $$  \endproclaim
\bigskip\flushpar
In this theorem the constants $C$ depend only on $n$, and integrals
and sums should be treated in the Cauchy sense. We will return to
functions of bounded variation while considering the radial case.

\head 3. Generalizations of the Bochner-Riesz means. \endhead
Let us consider another generalization of the linear means described
earlier. We mean the estimates (2.3) for more general
means of Bochner-Riesz type. They are strongly connected with a
support of the function generating these means.

{\bf 3.1.} The strongest estimates from above were obtained by Colzani and
Soardi \cite{CoS}. Their method is the direct generalization of that
used by Yudin \cite{Y1} and lies in the framework of our topic.

Suppose $S \subset \R$ is an open bounded set whose boundary $\ds$
has finite upper Minkowski measure. Let us consider complex-valued
bounded functions $\fa$ on $\R$ satisfying the following assumptions:
$$\fa (x) =0 \quad \text{if $x$ does not belong to} \quad S ; \tag3.1$$
there exist an integer $m \ge 0$ and real numbers $\a > -{1 \over 2}$
and $\beta > -{3 \over 2}$ such that
$$ \fa \in C^{m+1} (S) ;  \tag3.2$$
$$|D^{\xi} \fa (x)| \le C \rho (x,\ds)^{\a} \quad \text{if} \quad
\xi_1 +...+ \xi_n =m \quad \text{and} \quad x \in S; \tag3.3$$
$$|D^{\xi} \fa (x)| \le C \rho (x,\ds)^{\beta} \quad \text{if} \quad
\xi_1 +...+ \xi_n =m+1 \quad \text{and} \quad x \in S . \tag3.4$$
If (3.1)--(3.4) are satisfied with $m \ge 1$, $\fa$ must also satisfy
       $$ \fa \in C^{m-1} (\R) .  \tag3.5$$
Since $\fa$ is supposed bounded, we may assume $\a \ge 0$ whenever
$m=0$. Let
$$\gamma = \min \biggl(1, \a + {1\over 2}, \beta + {3 \over 2}\biggr).\tag3.6$$
If $\beta =-{1 \over 2}$ and $\a \ge {1 \over 2}$ , let
$$ \fa \in C^{m+2} (S) ;  \tag3.7$$
$$|D^{\xi} \fa (x)| \le C \rho (x,\ds)^{-{3 \over 2}} \quad \text{if}\quad
\xi_1 +...+ \xi_n =m+2 \quad\text{and}\quad x \in S . \tag3.8$$
\proclaim{Theorem 3.1} {\rm ([CoS])} Let $S$ be as above and $\fa$
satisfies {\rm (3.1)--(3.5)}, and, in addition, {\rm (3.7), (3.8)} when
$\beta =-{1 \over 2}$ and $\a \ge {1 \over 2}$.
Let $p_c = {2n \over n+2(m+\gamma)}$. Then for all $N > 2$:

a) If $m+\gamma \le {n \over 2}$
$$\aligned & \Vert L_N^{\fa} \Vert_p \le C_p N^{{n \over 2} - (m+\gamma)}\\
& \Vert L_N^{\fa} \Vert_p \le C_p N^{n{p-1 \over p}} \ln^{{1 \over p}} N \\
& \Vert L_N^{\fa} \Vert_p \le C_p N^{n{p-1 \over p}} \endaligned
\qquad \aligned & \text{if} \quad 1\le p < p_c , \\ & \text{if} \quad
p = p_c , \\ & \text{if} \quad p_c < p \le 2 .  \endaligned$$

b) If $m+\gamma > {n \over 2}$
$$\Vert L_N^{\fa} \Vert_1 \le C . $$          \endproclaim

M. Vignati [V] generalized these results to the case of nonisotropic
metrics in $\R.$ V. A. Yudin [Y3] showed that these estimates cannot
be asymptotically improved for $N\to\infty$ in the class of sets
considered.

{\bf 3.2.} Special examination of general conditions for lower estimates was
begun in \cite{Y2}, where the lower bound $\ln^n\!N$ for the order of growth
of \Lc\ of "all reasonable" partial sums is established, namely, for those
generated by convex sets which may contain a certain ball inside. Such an
investigation was continued, as it was mentioned above, in \cite{CaS}, and
then in [L2, L3] (see Theorem 1.3). The recent result from [LRZ]
generalizes the left-hand inequality in (2.3) in the spirit of Theorem 1.3.

Let $S=\operatorname{supp} \fa$ be the support of a function $\fa(x)$, where
$S$ is not necessarily a compactum. In what following we will be interested
in functions $\fa (x) = \fa_{r,\alpha} (x)$, which are $r$-smooth inside $S$,
and may be represented in a certain neighborhood of $\ds$ as follows:
$$\fa_{r,\a} (x) = f(x) (\rho (x))^{\a} , \tag3.9$$
where $f \in C^r (\R)$ and does not vanish on $\ds$, while $\rho (x) =0$ if
$x \not \in S$, and $\rho (x) = \rho (x,\ds)$ if $x \in S$.
Notice, that $\rho(x)$ is a smooth function in a neighborhood
of $\ds$ when $x\in S$ (see e.g., \cite{Gi, Appendix B}).
It should be mentioned that in \cite{CoS} the following obvious
consequence of (3.1)--(3.5) is proved:
\proclaim{Lemma 3.1} Suppose $S \subset \R$ is a bounded open set
such that $S$ has finite upper Minkowski measure and $\lambda$ is a
bounded complex-valued function on $\R$ satisfying {\rm (3.1)-(3.5).} Then
there exists a constant $C > 0$ such that
$$|\fa (x)| \le A \rho (x,\ds)^{\alpha + m} \ \ for \ \ all \ \ x \in
S . $$     \endproclaim

The following theorem shows that the generalizations cited
of the norms of the Bochner-Riesz means for upper
estimates and for lower estimates, respectively, are not far
one from another.
\proclaim{Theorem 3.2} {\rm ([LRZ])} Suppose that there exist an open set $U$
and a hypersurface $V$ of smoothness $r > \max ( 1, \Cp + \a )$, where
$0 \le \a < \Cp$, with non-vanishing principal curvatures, such that
$\ds \cap U = V$. Suppose, further, that in $U \cap S$ we have
$\fa (x) = \fa_{r,\a} (x)$.
 Then there exists a positive constant $C_{S,\fa}$ depending
only on $S$ and $\fa$ such that
$$\Vert L_{NS}^{\fa} \Vert \ge C_{S,\fa} N^{\Cp - \a} $$
for large $N$.              \endproclaim

{\bf 3.3.} We want to outline three knot points on which the proof of Theorem 7
is based. The first one was proved in a discussion of the author
and Belinskii.
\proclaim{Lemma 3.2} {\rm ([L2, LRZ])} Let $K$ be a set in $\R$ and
$\psi$ be a bounded measurable function with support in $K$. Then for
every point $x_0 \in \R$, for every ball $B_{\delta} (x_0)$ of radius
$\delta$ centred at $x_0$, and for every function $\ffi$ supported
in $B_{\delta} (x_0)$ and having the \Ft\ integrable over all $\R$,
there exists a constant $C$, depending only on $\ffi$, such that
$$\Vert L_K^{\psi} \Vert \ge C \Vert L_{K \cap B_{\delta} (x_0)}^{\psi
\ffi} \Vert . $$              \endproclaim
\demo{Proof}
We have
$$\Vert L_K^\psi \Vert =\operatorname{sup} \limits_{\Vert f \Vert
\le 1} \Vert L_K^\psi (f;\cdot)\Vert \ge \operatorname{sup}
\limits_{\Vert T_{B_{\delta}(x_0)}\Vert \le 1} \Vert L_K^\psi
(T_{B_{\delta}(x_0)};\cdot)\Vert, \tag3.10 $$
where $T_{B_{\delta}(x_0)}$ denotes all the trigonometric polynomials
with spectrum in $B_{\delta}(x_0)$.
According to \cite{Be1, Corollary 2} the following inequality holds
for every $f \in C(\T)$:
$$\left\Vert L_{B_{\delta}(x_0)}^\phi(f;\cdot) \right\Vert \le (2\pi)^{-n}
{\Vert \hat \phi \Vert}_{L_1(\R)} \Vert f \Vert .$$
Since the image of
$L_{B_{\delta}(x_0)}^\phi$ is only a part of all polynomials
$T_{B_{\delta}(x_0)}$, it follows from (3.10) that with some  constants
$$\align\Vert L_K^\psi \Vert \ge &\sup\limits_{\left\Vert
L_{B_{\delta}(x_0)}^\phi(f;\cdot) \right \Vert \le 1} \left\Vert L_K^\psi
\left(L_{B_{\delta}(x_0)}^\phi(f;\cdot);\cdot\right) \right\Vert \ge \\
&C \sup\limits_{\Vert f \Vert \le 1} \left\Vert
L_{K \bigcap B_{\delta}(x_0)}^{\psi \phi}(f;\cdot) \right\Vert=
C\left\Vert L_{K \bigcap B_{\delta}(x_0)}^{\psi \phi} \right\Vert .
\endalign$$
The lemma is proved.\qed     \enddemo

This lemma is of certain interest by itself, but mainly as a tool
for some estimates from below. It is more or less clear that after
its application we can pass to estimates (from below) of the \Lc\
of a method of summability, generated by a function with the support
possessing in global the properties which are local for the given set
in Theorem 7. A similar way to make ``global from local'' may be found in [Se].

The next step of the proof is the application
of Theorem 2.1, more precisely the lower estimate for $p=1$.
After that we need  appropriate asymptotic estimates of the \Ft\ of
the functions considered. Let us formulate in full such a
result strongly based on estimates of singularities of the Radon
transform, due to Ramm and Zaslavsky (see [RZ1, RZ2]).
\proclaim{Theorem 3.3} {\rm([LRZ], see also [RZ] and [RK])} Let $S$ be the
compact support of a function $\fa (x) = \fa_{r,\alpha} (x)$ with
$\alpha \ge 0$ and $r > \max (1,\Cp + \alpha)$. Let $S$ be convex,
with the $r$-smooth boundary $\ds$, and suppose the principal curvatures
of $\ds$ never vanish. Let $\theta \in \R$ be a vector on the unit
sphere, $x^+ (\theta)$ and $x^- (\theta)$ be the (uniquely defined)
points of $\ds$ at which the function $\theta_1 x_1 +...+ \theta_n
x_n $ attains maximum and minimum on $\ds$, respectively. Then
for $t \to +\infty$
$$\faf (t\theta) = t^{-\a - {n+1 \over 2}} \left(
\Xi^+ (\theta) e^{itx^+ (\theta)\theta} + \Xi^- (\theta)
e^{itx^- (\theta)\theta} + o(1) \right) ,$$
$$\Xi^{\pm} (\theta) = (2\pi)^{\Cp} \Gamma (\a + 1) e^{\pm i\pi
{2\a +n+1 \over 4}} f(x^{\pm} (\theta)) (\varkappa^{\pm}
(\theta))^{-{1 \over 2}} , $$
where the remainder term is small uniformly in $\theta,$ and
$\varkappa^{\pm} (\theta)$ is the Gaussian curvature of $\ds$
at the points $x^{\pm} (\theta)$, respectively.
\endproclaim
This result continues and develops the well-known asymptotic estimate for
the characteristic function of a convex set \cite{GGV}. There is an "almost
all" gap between Theorem 8 and the result in \cite{P6} in the
two-dimensional case. We must mention that many authors use one result
of Herz \cite{Hz} for estimates of Fourier transforms. But smoothness
assumptions in this work are essentially more restrictive than those in
\cite{GGV} (and of course in Theorem 3.3) since the author was interested
in sharp estimate for the remainder term. This explains, for example,
the smoothness conditions in \cite{CaS} or [Br1, Br2].
Of course we omit a number of technical details and estimates,
not easy by no means.

{\bf 3.4.} Let us pay our attention to the following circumstance. One can see
that in many results cited the value $\Cp$ behaved like a star actor on the
stage. It is not an occasional event, and this number is called "critical
order" for the Bochner-Riesz means. Let us compare Theorem 2.2 with the
following well-known result of Stein:
\proclaim{Theorem 3.4} {\rm ([S1])} The following asymptotic formula holds:
$$\Vert L_N^{R_{\Cp}} \Vert_{L_1 (\T) \to L_1 (\T)} =
\omega_n \ln N + O( 1 ) .  \tag3.11 $$      \endproclaim
\bigskip\flushpar
This asymptotics was obtained as a corollary to some general estimates
of the difference between the corresponding kernel
$$\sum \limits_{|k| \le N} R_{\Cp} \biggl({k \over N}\biggr)\, \ek$$
and its  integral analog. The constant $\omega_n$ was not indicated explicitly.
Here the Lebesgue constants of the Bochner-Riesz means lose their
power rate of growth, and behave as the \Lc\ of one-dimensional partial
sums. This likeness is not casual.
Before formulating one recent generalization of Theorem 3.4 we want
to mention that (3.11) is a relatively simple corollary to Theorem 2.1.
\proclaim{Corollary 2.2} The following asymptotic formula holds:
$$||L_N^{R_{\Cp}}||_{L_1 (\T)
\to L_1 (\T)} =\frac{4\Gamma ({n+1 \over 2})}{\pi^{{3 \over 2}}
\Gamma ({n \over 2})}\ln N + O( 1 ) .$$   \endproclaim
\demo{Proof} {\rm ([Be2])} By Theorem 2.1 and (2.5), (2.6), and (2.8) we have
$$||L_N^{R_{{n-1\over 2}}}||=\pi^{-{n+1\over 2}}\Gamma({n+1\over 2})
\int\limits_{1\le|x|\le N}\biggl|\prod\limits_{j=1}^n{x_j\over
2N\sin{x_j\over 2N}}{\cos(|x|-{\pi n\over 2})\over |x|^n}\biggr|\,dx+O(1).$$
The following relation $${x_j\over 2N\sin{x_j\over 2N}}-1=O\biggl(\biggl|
{x_j\over N}\biggr|^2\biggr)$$ and estimates from the proof of
Corollary 2.1 imply
$$||L_N^{R_{{n-1\over 2}}}||=\pi^{-{n+1\over 2}}\Gamma({n+1\over 2})
\int\limits_{1\le|x|\le N}\biggl|
{\cos(|x|-{\pi n\over 2})\over |x|^n}\biggr|\,dx+O(1).$$
We obtain after passage to spherical coordinates
$$\align ||L_N^{R_{{n-1\over 2}}}||&=\pi^{-{n+1\over 2}}\Gamma({n+1\over 2})
\int\limits_{1\le|x|\le N}\biggl|
{\cos(|x|-{\pi n\over 2})\over |x|^n}\biggr|\,dx+O(1)\\ \quad\\
&=\pi^{-{n+1\over 2}}\Gamma({n+1\over 2}){2\pi^{{n\over 2}}\over
\Gamma({n\over 2})}\int\limits_1^N \biggl|
{\cos(t-{\pi n\over 2})\over t}\biggr|\,dt+O(1).\endalign$$
It is well-known that the last integral is $${2\over\pi}\ln N+O(1)$$
(see e.g., [Z], Vol.1, Ch.2), and this completes the proof.  \qquad\qed\enddemo

Furthermore, again Theorem 2.1 and certain technique connected with
Theorem 3.3 allow us to obtain the mentioned generalization as follows.
\proclaim{Theorem 3.5} {\rm([L5])} Let $S$ be the compact support
of a function $\lambda=\lambda_{n,{n-1\over 2}},$ with the
$n$-smooth boundary $\ds$.
Assume that $S$ is convex and the
principal curvatures of $\ds$ never vanish. Then there exists a positive
constant $C_{S,\fa}$ depending only on $S$ and $\fa$ such that
$$ \Vert L_N^{\fa} \Vert_{L_1 (\T) \to L_1 (\T)} =
C_{S,\fa} \ln N + o( \ln N ) \tag3.12 $$
for large $N$.                  \endproclaim
\remark{Remark 3.1} The following formula is given in [L5] to calculate
$C_{S,\lambda}:$
$$C_{S,\lambda}=(2\pi)^{{n+3\over 2}}\Gamma({n+1\over 2})\int\limits_
{|\theta|=1}\,d\theta\int\limits_0^{2\pi}|(-1)^n\phi^+(\theta)e^{it}
+\phi^-(\theta)|\,dt$$ where $\phi^{\pm}(\theta)=f(x^{\pm}(\theta))
(\varkappa^{\pm}(\theta))^{-{1\over 2}}$ (cf. Theorem 3.3). Simple
calculations yield for the Lebesgue constants of the usual Bochner-Riesz
means the same constant as in Corollary 2.2.    \endremark
\remark{Remark 3.2} It is obvious that if, in the conditions of Theorem 10,
to take $\fa = \fa_{r,\a}$, with $r > n$ and $\a > \Cp$, we will obtain
$\Vert L_N^{\fa} \Vert = O(1)$ (cf. b) in Theorem 5).  \endremark

\head 4. "Radial" results. \endhead
A lot of attention to the Bochner-Riesz means and certain of their
generalizations has been given. But we have not speak yet about one more
peculiarity of the Bochner-Riesz means. We mean the fact that they
are generated by the function $R_{\a}$ which is radial, that is
depending only on $|x|$. Such functions play a special role in Fourier
Analysis, and there are many ways to exploit the radiality.

{\bf 4.1.} Now we want to consider one special class of radial  functions,
close, in many respects, to the Bochner-Riesz means of critical order.

Let $\fa (x) = \fa_0 (|x|)$ be a radial function satisfying the following
conditions:
$$\fa_0 \in C [0,\infty), \qquad \fa_0 \in C^{[{n-2 \over 2}]}
(0,\infty) , \tag4.1$$
$$t^r \fa_0^{(r)} (t) \to 0 \ \ \text{as} \ \ t \to 0, \quad r=1,2,...,
[{n-2 \over 2}] \quad ( n > 3 ) ,  \tag4.2$$
$$\Lambda (t) = t^{\Cp} \fa_0^{(\Cp)} (t), \quad
t^r \fa_0^{(r)} (t) \to 0 \ \ \text{as} \ \ t \to \infty , \quad r=0,1,2,...,
[{n-2 \over 2}] , \tag4.3 $$
$$\Lambda \ \ \text{has a bounded variation} \ \ V_{\Lambda} \ \ \text{on} \ \
[0,\infty) . \tag4.4 $$

Here the fractional derivative is understood in the Weil sense (see e.g.,
\cite{BE}), namely, for a function $g$ defined on $[0,\infty)$
$$W_{\a} (g;t) = {1 \over \Gamma (\a)} \int\limits_t^\infty  g(u)
(u-t)^{\a -1} \, du $$
is the Weil integral of fractional order $\a$, and for $0 < \a < 1$
$$g^{(\a)} (t) = {d \over dt} W_{1-\a} (g;t) $$
is the fractional derivative of order $\a$. For $\gamma = r+\a$,
$r=1,2,...$, $0 <\a < 1$,
$$g^{(\gamma)} (t) = {d^r \over dt^r} g^{(\a)} (t) $$
is the fractional derivative of order $\gamma$.
\bigskip
\proclaim{Theorem 4.1} {\rm ([BL1, BL2])} Let $\fa (x) = \fa_0 (|x|)$
be a radial function satisfying {\rm (4.1)--(4.4).} Then
$$\Vert L_N^{\fa} \Vert_{L^1 (\T) \to L^1 (\T)} = (2\pi)^{-n}
\int\limits_{|x| \le \pi N} |\faf (x)| \, dx + O ( V_{\Lambda} +
|\fa (0)| ) .  \tag4.5$$        \endproclaim
\demo{Proof} It suffices to prove that the series
$$\sum\limits_k \fa\biggl({k\over N}\biggr)e^{ikx} \tag0.3$$
is the Fourier series of an integrable function in order
(0.4) to be true. Consider
$$\align &\biggl| ||L_N^{\fa}||-(2\pi)^{-n}\int\limits_{|x|\le\pi N}
|\faf(x)|\,dx-(2\pi)^{-n}|\fa(0)|\biggr| \\
&=R_N\le(2\pi)^{-n}\int\limits_{\T}\biggl|\sum\limits_{k\ne 0}\fa\biggl(
{k\over N}\biggr)e^{ikx}-\Phi_N(x)\biggr|\,dx, \endalign$$   where
$$\Phi_N(x)=\cases N^n\faf(Nx), & |x|\le\pi,\\ 0, & x\in\T\setminus
\{x: |x|\le\pi\}. \endcases$$ \medskip\flushpar
Let us calculate the $k$-th Fourier coefficient of the periodic in
each variable continuation of this function. No confusion will result
saving the same notation. In [BL1, BL2] the following results
were obtained: if a function $\fa$ satisfies conditions (4.1)--(4.4)
its Fourier transform can be calculated by the following relation:
$$\faf(x)={(2\pi)^{{n\over 2}}(-1)^{[{n\over 2}]}\over\Gamma({n-1\over 2})}
|x|^{1-{n\over 2}}\int\limits_0^\infty\Lambda(t)t^{{n\over 2}}Q(|x|t)\,dt,
\tag4.6$$ where $Q(r)=\int\limits_0^1(1-s)^{{n-3\over 2}}s^{{n\over 2}}
J_{{n\over 2}-1}(rs)\,ds,$ and the inverse formula holds:
$$\fa(x)=\lim\limits_{A\to\infty}(2\pi)^{-n}\int\limits_{|u|\le A}
\faf(u)e^{ixu}\,du.\tag4.7$$ For generalizations of this result
see [L6], [L7]. For $|k|>0$ (4.6) and (4.7) yield
$$\align\widehat\Phi_N(k)&=(2\pi)^{-n}\int\limits_{|u|\le\pi}
N^n\faf(Nu)e^{-iku}\,du \\
&=(2\pi)^{-n}\int\limits_{|u|\le\pi N}\faf(u)e^{-i{k\over N}u}\,du\\
&=\fa({-k\over N})-(2\pi)^{-n}\int\limits_{|u|>\pi N}\faf(u)
e^{-i{k\over N}u}\,du.\endalign$$
For $k=0,$ the passage to spherical coordinates and (4.6) yield
$$\widehat\Phi_N(0)=(2\pi)^{-n}\int\limits_{|u|\le\pi N}\faf(u)\,du=
{2^{n-1}(-1)^{[{n\over 2}]}\over\Gamma({n\over 2})\Gamma({n-1\over 2})}
\int\limits_0^{\pi N}r^{{n\over 2}}\,dr\int\limits_0^\infty\Lambda(t)
t^{{n\over 2}}Q(rt)\,dt.$$
Use the well-known formula (see e.g., [BE], 7.2.8(50),(51))
$${d\over dt}t^{\pm\nu}J_{\nu}(t)=\pm t^{\pm\nu}J{\nu\mp1}(t)\tag4.8$$
and denote $q(r)=\int\limits_0^1(1-s)^{{n-3\over 2}}s^{{n\over 2}-1}
J_{{n\over 2}}(rs)\,ds.$ This implies
$$ \widehat\Phi_N(0) ={2^{n-1}(-1)^{[{n\over 2}]}\over
\Gamma({n\over 2})\Gamma({n-1\over 2})}(\pi N)^{{n\over 2}}
\int\limits_0^{\pi N}r^{{n\over 2}}\,dr\int\limits_0^\infty\Lambda(t)
t^{{n\over 2}-1}\biggl[q(\pi Nt)-\a_2(\pi Nt)^{-{n\over 2}}
\biggr].$$             For $q$ the following asymptotic relation
was obtained in [BL2] (see also [L6, L7]):
$$q(r)=\a_1 r^{-{n-1\over 2}}J_{n-{1\over 2}}(r)+\a_2 r^{-{n\over 2}}
+O(r^{-{n+2\over 2}})\tag4.9$$
as $r\to\infty,$ where $\a_1$ and $\a_2$ are some constants.
Integrating by parts one obtains
$$\align \widehat\Phi_N(0) &={2^{n-1}(-1)^{[{n\over 2}]}\over
\Gamma({n\over 2})\Gamma({n-1\over 2})}(\pi N)^{{n\over 2}}
\left\{\Lambda(t)\int\limits_0^t r^{{n\over 2}-1}
\biggl[q(\pi Nr)-\a_2(\pi Nr)^{-{n\over 2}}\biggr]\,dr\right|_0^\infty\\
\quad\\&-\left.\int\limits_0^\infty\biggl\{\int\limits_0^t r^{{n\over 2}-1}
\biggl[q(\pi Nr)-\a_2(\pi Nr)^{-{n\over 2}}\biggr]\,dr\biggr\}\,d\Lambda(t)
\right\}.\endalign$$
In order to get $|\widehat\Phi_N(0)|\le CV_{\Lambda}$ it suffices,
taking into account (4.4), to prove the boundedness of the value
$$\align \sup\limits_{N,t}(\pi N)^{{n\over 2}}&\biggl|\int\limits_0^t
r^{{n\over 2}-1}\biggl[q(\pi Nr)-\a_2(\pi Nr)^{-{n\over 2}}\biggr]\,dr\biggr|\\
\quad\\ = \sup\limits_{N,t:\pi Nt>1}&\biggl|\int\limits_1^{\pi Nt}
r^{{n\over 2}-1}\biggl[q(\pi Nr)-\a_2(\pi Nr)^{-{n\over 2}}
\biggr]\,dr\biggr|+O(1).\endalign$$
It follows from (4.9) that the right-hand side is equal to
$$\align \sup\limits_{N,t:\pi Nt>1}\ &\biggl|\a_1\int\limits_1^{\pi Nt}
r^{-{1\over 2}}J_{n-{1\over 2}}(r)\,dr+O(
\int\limits_1^{\pi Nt}r^{-2}dr)\biggr|+O(1) \\ \quad\\
=\sup\limits_{N,t:\pi Nt>1}\ &\biggl|\a_1\int\limits_1^{\pi Nt}
r^{-{1\over 2}}J_{n-{1\over 2}}(r)\,dr\biggr|+O(1).\endalign$$
The asymptotic formula (2.6) makes the claimed estimate obvious. Thus
$$R_N\le(2\pi)^{-n}\int\limits_{\T}\biggl|\sum\limits_{k\ne 0}\biggl[
\quad\int\limits_{|u|>\pi N}\faf(u)
e^{-i{k\over N}u}\,du\biggr]e^{ikx}\biggr|\,dx+CV_{\Lambda}.$$
Apply the Cauchy-Schwarz inequality to the outer integral and
the Parseval equality. We get
$$R_N\le C\biggl\{\sum\limits_{k\ne 0}\biggl|\ \int\limits_{|u|>\pi N}
\faf(u)e^{-i{k\over N}u}\,du\biggr|^2\biggr\}^{{1\over 2}}+CV_{\Lambda}.$$
The Cauchy-Poisson formula (see e.g.,[Bc1], Th.56) and (4.6) yield
$$R_N\le CN^{{n\over 2}-1}\biggl\{\sum\limits_{k\ne 0}|k|^{2-n}
\biggl|\int\limits_{\pi N}^\infty J_{{n\over 2}-1}\biggl({|k|\over N}r\biggr)
r\,dr \int\limits_0^\infty \Lambda(t)t^{{n\over 2}}Q(rt)\,dt\biggr|^2
\biggr\}^{{1\over 2}}+CV_{\Lambda}.$$
Integration by parts in $t$ implies:
$$\align R_N & \le CN^{{n\over 2}-1}\left\{\sum\limits_{k\ne 0}|k|^{2-n}
\biggm|\int\limits_{\pi N}^\infty J_{{n\over 2}-1}
\biggl({|k|\over N}r\biggr)\,dr\Lambda(t)t^{{n\over 2}}q(rt)\right|_0^\infty \\
&-\left.\int\limits_{\pi N}^\infty J_{{n\over 2}-1}\biggl({|k|\over N}r\biggr)
\,dr \int\limits_0^\infty t^{{n\over 2}}q(rt)\,d\Lambda(t)\biggm|^2
\right\}^{{1\over 2}}+CV_{\Lambda}.\endalign$$
After applying generalized Minkowski's inequality and (4.4)
we get, as above, that it suffices to prove the boundedness of the value
$$\align \sup\limits_{N,t}N^{{n\over 2}-1}&\left\{\sum\limits_{k\ne 0}
|k|^{2-n}\biggl|t^{{n\over 2}}\int\limits_{\pi N}^\infty J_{{n\over 2}-1}
\biggl({|k|\over N}r\biggr)q(rt)\,dr\biggr|^2\right\}^{{1\over 2}} \\ \quad\\
=\sup\limits_{N,t}N^{{n\over 2}-1}&\left\{\sum\limits_{k\ne 0}|k|^{2-n}
\biggl|t^{{n\over 2}}\int\limits_{\pi N}^\infty r^{{n\over 2}}
J_{{n\over 2}-1}\biggl({|k|\over N}r\biggr)\,dr\int\limits_0^1(1-s)^{{n-3\over
 2}} s^{{n\over 2}-1}{J_{{n\over 2}}(rts)\over r^{{n\over 2}}}ds
\biggr|^2\right\}^{{1\over 2}}.\endalign$$
Integrate by parts using (4.8) and obtain
$$\align &\sup\limits_{N,t}N^{{n\over 2}} \biggl\{\sum\limits_{k\ne 0}
|k|^{-n} \biggm| t^{{n\over 2}} J_{{n\over 2}}
\biggl({|k|\over N}r\biggr)q(rt) \biggr|_{\pi N}^\infty \\ \quad\\
&+ t^{{n\over 2}+1}\int\limits_{\pi N}^\infty  J_{{n\over 2}}\biggl
({|k|\over N}r\biggr)\,dr\int\limits_0^1(1-s)^{{n-3\over 2}}
s^{{n\over 2}}J_{{n\over 2}+1}(rts)ds
\biggl|^2\biggr\}^{{1\over 2}}.\endalign$$
Relations (4.8) and (4.9) and convergence of the series $\sum\limits_
{k\ne0}|k|^{-n-1}$ imply the boundedness of the integrated terms. Further,
integration by parts and (4.10) yield
$$\int\limits_0^1(1-s)^{{n-3\over 2}}s^{{n\over 2}}J_{{n\over 2}+1}(rts)ds
=\a_3 r^{-{n\over 2}}\sin(r-{\pi n\over 2})+O(r^{-{n+2\over 2}}).$$
Estimates with the remainder term are trivial. Let us estimate
$$\sup\limits_{N,t}N^{{n\over 2}}\left\{\sum\limits_{k\ne 0}
|k|^{-n}\left|t\int\limits_{\pi N}^\infty r^{-{n\over 2}}
J_{{n\over 2}}\biggl({|k|\over N}r\biggr)\sin(rt-{\pi n\over 2})\,dr
\right|^2\right\}^{{1\over 2}}.$$
Again integrate by parts and get
$$\align &\sup\limits_{N,t}N^{{n\over 2}}\left\{\sum\limits_{k\ne 0}
|k|^{-n}\left|r^{-{n\over 2}}J_{{n\over 2}}\biggl({|k|\over N}r\biggr)
\cos(rt-{\pi n\over 2})\right|_{\pi N}^\infty \right.\\ \quad\\
&-{|k|\over N}\int\limits_{\pi N}^\infty r^{-{n\over 2}}
J_{{n\over 2}+1}\biggl({|k|\over N}r\biggr)\cos(rt-{\pi n\over 2})\,dr
\biggl|^2\biggr\}^{{1\over 2}}.\endalign$$
The integrated terms are symply estimated. Apply now (2.6) to the
last integral. Estimates for the remainder term are obvious. Using
also simple trigonometric identities we have to estimate
$$\sup\limits_{N,t}N^{{n-1\over 2}}\left\{\sum\limits_{k\ne 0}
|k|^{1-n}\biggl|\ \int\limits_{\pi N}^\infty r^{-{n+1\over 2}}
\sin r({|k|\over N}-t)\,dr\biggr|^2\right\}^{{1\over 2}}.\tag4.10$$
Notice that estimates for similar values with $\sin r({|k|\over N}+t)$
or $\cos r({|k|\over N}\pm t)$ on place of $\sin r({|k|\over N}-t)$
are the same. Assume that $Nt$ is big enough. Split the sum in (4.10)
into three ones: $1\le |k|<Nt-1,$ $Nt-1\le|k|\le Nt+1,$ and $Nt+1<|k|<\infty.$
Integration by parts implies for the integral in (4.10) the estimate
$$N^{-{n+1\over 2}}\biggl|{|k|\over N}-t\biggr|^{-1}=
N^{-{n-1\over 2}}\biggl|\,|k|-Nt\biggr|^{-1}.$$
Therefore the boundedness of the following sums:
$$\sum\limits_{1\le|k|<Nt-1}|k|^{1-n}(Nt-|k|)^{-2} \qquad\text{and}
\qquad\sum\limits_{Nt+1<|k|<\infty}|k|^{1-n}(|k|-Nt)^{-2}$$
has to be established when estimating with respect to the first and
third domains. This is easy to see passing to integrals instead of sums.
For the second one we obtain
$$\align\sup\limits_{N,t}\,N^{{n-1\over 2}}&\left\{\sum\limits_{Nt-1\le
|k|\le Nt+1}|k|^{1-n}\biggl|\int\limits_{\pi N}^\infty r^{-{n+1\over 2}}
\sin r({|k|\over N}-t)\,dr\biggr|^2\right\}^{{1\over 2}}\\
\le  N^{{n-1\over 2}}&\left\{\sum\limits_{Nt-1\le|k|\le Nt+1}|k|^{1-n}
\biggl|\int\limits_{\pi N}^\infty r^{-{n+1\over 2}}\,dr
\biggr|^2\right\}^{{1\over 2}}\\
\le  C &\left\{\sum\limits_{Nt-1\le|k|\le Nt+1}|k|^{1-n}
\right\}^{{1\over 2}} \le C. \endalign$$
When $Nt$ is small similar estimates valid after splitting the sum
into two ones: $1\le|k|\le 3$ and $3<|k|<\infty.$ The proof is
complete. \qquad\qed\enddemo
\remark{Remark 4.1} One can see that (0.3) is the Fourier series of a function
not only from $L^1(\T)$ but from $L^2(\T)$ as well. \endremark
\bigskip \remark{Remark 4.2}
Observe that besides other applications, say, to approximation
on the class of functions with bounded polyharmonic operator,
Theorem 4.1 allows to obtain (3.11) as a simple corollary once more.
Indeed, conditions (4.1)--(4.4) are verified easily. Then the
estimates are similar to those in the proof of Corollary 2.2, and of
course with the same constant.       \endremark

So we have several different approaches with (3.11) as intersection.
This fact cannot be senseless?!

{\bf 4.2.} One of the features of radial functions is that they combine, in
a certain sense, some properties of the multi-dimensional case and some
properties of the one-dimensional case. In our situation it may be expressed
by the following relation.
\proclaim{Theorem 4.2} {\rm ([BL1])} Let $\fa$ be a function satisfying
{\rm (4.1)--(4.4)} and, moreover,
$$\int\limits_0^1 {|\Lambda (t)| \over t} \, dt < \infty . \tag4.12$$
Then we have for $n \ge 2$
$$\align \int\limits_{1 \le |x| \le N} |\faf (x)| \, dx & = \frac
{2^{{n+3 \over 2}} \pi^n}{\Gamma ({n \over 2})} \int\limits_1^N \left|
\int\limits_0^{\infty} \Lambda (t) \sin (st -{\pi n \over 2}) \, dt
\right| \, ds \\
& + O ( V_{\Lambda} + ||\Lambda ||_{C[0,\infty)} + \int\limits_0^1
{|\Lambda (t)| \over t} \, dt ) . \tag4.13  \endalign$$
\endproclaim

\remark{Remark 4.3} The condition (4.12) is sharp and cannot be
removed.\endremark
\flushpar
Hence, we can apply many one-dimensional results dealing with the behavior
of \Ft s to such functions $\fa$. Let us give such an example.
\proclaim{Proposition 4.1} {\rm ([BL1])} If $\fa$ satisfy {\rm (4.1)-(4.4)}
and {\rm (4.12)}, and $\Lambda$ has at least one point of discontinuity, then
$$\Vert L_N^{\fa} \Vert_{L^1 (\T) \to L^1 (\T)} = M(\Lambda)
\ln N+o(\ln N) ,$$ where $M(\Lambda)$ is an average of some almost
periodic function built in accordance with $\Lambda$.   \endproclaim

This result as well as Theorem 4.1 are generalizations of one-dimensional
results in \cite{Be0}.

Some other conditions allowing a similar connection between the radial case
and the one-dimensional one were obtained by Podkorytov. Let us formulate them.
\proclaim{Theorem 4.3} {\rm ([P4])} Let $\fa_0 \in C[0,\infty)$
and $\fa_0 (t) = 0$ for $t \ge 1$. Then the following integrals
converge simultaneously:
$$ \int\limits_{\R} |\faf (x)| \, dx $$
and
$$\int\limits_0^{\infty} s^{\Cp} \left| \int\limits_0^1 t^{\Cp} \fa_0 (t)
\cos (2\pi st - {\pi (n-1) \over 4}) \, dt \right| \, ds . $$ \endproclaim
\bigskip
If to compare the latter two theorems, one can see that conditions
(4.1)--(4.4) and (4.12) is the price payed for $\fa$ should not be
necessarily boundedly supported and having the integrable \Ft\.

{\bf 4.3.} Let us describe one more "radial" result due to Trigub. The estimate
(2.11) plays an important role in its proof. Consider a function $\fa_0 (t)$
in $[0,\pi]$ and expand it in the cosine series:
$$\fa_0 (t) \sim \sum_{j=0}^{\infty} a_j \cos jx .$$
\proclaim{Theorem 4.4} {\rm ([T2])} Let $\fa_0 \in C^{[\Cp]} [0,\pi]$,
and $\fa_0^{(r)} (\pi) = 0$, $0 \le r \le [\Cp]$. Then
$$\sup_N \int \limits_{\T} \biggl| \sum_{|k| \le N} \fa_0 \biggl(
{|k|\pi \over N}\biggr) \ek
\biggr| \, dx \le C \sum_{j=0}^{\infty} j^{\Cp} |a_j| \ln (j+1) . \tag4.14$$
\endproclaim
\bigskip
It is supposed that the series on the right-hand side converges and
$C$ depends only on $n$. It suffices to supplement (4.14) with the
condition $\fa_0 (0) = 1$ for the summability  on the whole class of
periodic continuous functions.

It may be shown that for $a_j$, $j \ge 1$, with alternating signs,
the opposite inequality holds provided $n = 1 (\operatorname{mod}4)$.
This generalizes the corresponding one-dimensional result \cite{T1}.
It is possible to consider coefficients $b_j$ of the sine expansion
of $\fa_0$, instead of $a_j$, as well.

\head 5. "Polyhedral" results. \endhead
{\bf 5.1.} Let us quote another extract from the book [DC]: "Concentric
polygons are an obvious thing to try, but this turns out to be no more
interesting than repeating several one-dimensional results. It doesn't give
any new mathematics, and it avoids having to think deeply about Fefferman's
result. \footnote{The famous solution of the multiplier problem for the ball
in \cite{F2} is meant.} To avoid thinking about a subject is almost always
a mistake; at best you are in for some big surprises later on".

This passage is not absolute truth even if one speaks about polygons (or,
maybe more exactly, polyhedra) with the sides parallel to coordinate planes.
In this case one frequently obtains nothing more than the product of
one-dimensional estimates. But even there the above-mentioned quotation
contradicts the fact that there exists another Fefferman's bright result
\cite{F1}, which gives an example of continuous function with rectangularly
divergent partial sums.  And if to consider more general objects
being within the frame of "polyhedral" case, one can find many non-trivial
problems. We will touch some of them connected with our topic.

We must say that, in general, this case has a "logarithmic" nature.
For the case of arbitrary parallelogram it was shown by Belinskii \cite{Be2};
that the same method proves the more general
case was also mentioned here. This was realized later in [P3, Bb].
More precisely, there exist two positive constants $C_1$, $C_2$,
with $C_1 < C_2$, such that for each polyhedron $E$
$$C_1 \ln^n N \le \int\limits_{\T} | \sum_{k \in N\!E} \ek | \, dx
\le C_2 \ln^n N . \tag5.1$$
Thus, we see an essential difference between this case and the "ball"
case (1.2). In the latter case the Lebesgue constants have a power growth.
We intend to pay our attention to two important problems around polyhedra.

{\bf 5.2.} One of them deals with a quite natural question: whether some sets
of spectrum of partial sums exist for which the norms of the corresponding
operators (1.1) have an intermediate (between (5.1) and (1.2)) rate of
growth with respect to $N$-dilations of these sets. The positive answer
to this question was given by Podkorytov (similar results were given
in [YY2]).

Let $C_1$ and $C_2$ denote, as above, positive constants such that $C_1 < C_2$.
\proclaim{Theorem 5.1} {\rm ([P4])}
{\bf 1.} For any $p>2$ there exists a compact, convex set $E$ for which
$$C_1 \ln^p N \le \int\limits_{\bold T^2} | \sum_{k \in N\!E} \ek | \, dx
\le C_2 \ln^p N , \qquad N \ge 2 .   \tag5.2$$
\flushpar
{\bf 2.} For any $p \in (0,{1 \over 2})$ and $\a > 1$ there exists a compact,
convex set $E$ for which
$$C_1 N^p \ln^{-\a p} N \le \int\limits_{\bold T^2} | \sum_{k \in N\!E}
\ek | \, dx \le C_2 N^p \ln^{2-2p} N , \qquad N \ge 2 .   \tag5.3$$
\endproclaim

It is well worth noticing that these sets look like polyhedra with an
infinite number of sides (of course, becoming smaller and smaller),
and the sketch of the proof, outlined when proving Theorem 1, was realized
there with taking into account very subtle technical peculiarities.

{\bf 5.3.} The second question is rather natural too and asks whether it is
possible to write a certain asymptotic relation instead of (5.1). Some partial
 cases were investigated by Daugavet \cite{D}, Kuznecova \cite{Ku}, Skopina
[Sk0, Sk2]. For example Kuznetsova generalized Daugavet's
result as follows. What tells both these results from many others is that not
dilations of certain fixed domain are taken. Namely, the following
result is true.
\proclaim{Theorem 5.2} {\rm ([Ku1, Ku2])} Let $B_{N_1,N_2}=\{(k_1,k_2):
{|k_1|\over N_1}+{|k_2|\over N_2}\le 1\}.$ The asymptotic equality
$$||S_{B_{N_1,N_2}}||=32\pi^{-4}\ln N_1\ln N_2 -16\pi^{-4}\ln^2 N_1
+O(\ln N_2)$$ holds uniformly with respect to all natural $N_1,$
$N_2,$ and $l={N_2\over N_1}.$      \endproclaim
\medskip\flushpar The case $l=1$ is the mentioned result of Daugavet.

An unexpected result was obtained again by
Podkorytov \cite{P6}. He has shown that there are two main cases. The first
one, to which the afore-mentioned asymptotic results may be referred,
deals with the polygons (we are speaking about two-dimensional results)
with integer, or rational vertices. Of course, those times any constant, the
same for all vertices, are appropriate. In this case one can show that the
estimates change insignificantly if, instead of sums, to consider the
corresponding integrals, that is, the \Ft\ of the indicator function
of the corresponding set. This circumstance allows to obtain the logarithmic
asymptotics, namely, $\int \limits_{\bold T^2} |\hat \chi_{N\!E} (x)|dx$
is equivalent to $C \ln^2 N.$

In the second case, that is, for any dilation at least one vertex has some
irrational coordinates, the situation changes qualitatively: the upper
limit and the lower limit, as $N \to \infty$, are different, and a limit
of $\int_{\T} |\sum_{k \in N\!E} \ek | \, dx $, as $N \to \infty$, does not
exist. In other words, in this case the behavior of the \Ft\ of the
indicator function of $N\!E$ does not express the nature of the
behavior of the corresponding partial sums. In \cite{P6} the quantitative
estimate of this phenomenon is given at once. Namely, for the triangles
$$E = E_{\a} = \{ (u,v):\qquad 0 \le u \le 1,\qquad 0 \le v \le \a u \}$$
the following theorem holds.
\proclaim{Theorem 5.3}
\flushpar
{\bf 1.} $\int \limits_{\bold T^2} \biggl| \sum \limits_{k \in N\!E_{\a}} \ek
\biggr| \, dx = \int \limits_{\bold T^2} | \hat \chi_{N\!E_{\a}} (x)| \, dx
+ \int \limits_0^{2\pi} \biggl| \sum \limits_{j=0}^N \{ \a j\}
e^{ijt} \biggr| \, dt + O(\ln N \ln \ln N)$
\medskip
\qquad where $\{...\}$ is the symbol of the fractional part.
\medskip
{\bf 2.} There exists irrational $\a$ such that
$$ \varlimsup \limits_{N \to \infty} {1 \over \ln^2 N} \int\limits_0^{2\pi}
\biggl| \sum \limits_{0 \le j \le N}\{ \a j\}e^{ijt}\biggr|\, dt> 0.\tag5.4 $$
\endproclaim
The main defect of this theorem is that it is true only for $\a$
from very scarce setand nothing know about other $\a.$ In a
recent paper of F. Nazarov and A. Podkorytov [NP] this uncertainty
is partly removed. Namely, the following is true. Denote by $I_N(\a)$
the integral from (5.4).
\proclaim{Statement 5.1} Let $\a$ be irrational. Then
\medskip
{\bf 1.}\qquad $0<C_1\le\varlimsup\limits_{N\to\infty}
{I_N(\a)\over\ln^2N}\le C_2.$

{\bf 2.} $\varliminf\limits_{N\to\infty}{I_N(\a)\over\ln^2N}=0$
if and only if$\a$ is a Liouville number, that is if and only if for
each $M>0$ there exist fractions ${p\over q}$ ($q\ge 2$) such that
$|\a-{p\over q}|\le{1\over q^M}.$

{\bf 3.} If $|\a-{p\over q}|\le{1\over q^M}$ for some $M>2$ and as
many fractions ${p\over q}$ as infinite ($q\ge 2$), then the fraction
${I_N(\a)\over \ln^2N}$ has no limit as $N\to\infty.$

{\bf 4.} The integral $I_N(\a)$ is concentrated on a set of small
measure, namely for all $N\ge 2$ and $\a$ irrational there exists
a set $E=E(N,\a)\subset\bold T$ such that $\operatorname{mes}(E)\le
e^{-\sqrt{\ln N}}$ while
$$\int\limits_{\bold T\setminus E} \biggl| \sum \limits_{0 \le j \le N}
\{ \a j\}e^{ijt}\biggr|\, dt\le C\ln^{{3\over 2}}N.$$

{\bf 5.} There exist numbers $0<\omega\le\Omega<\infty$ such that for
almost all $\a$
$$\omega=\varliminf\limits_{N\to\infty}{I_N(\a)\over\ln^2N}\qquad
\text{and}\qquad\Omega=\varlimsup\limits_{N\to\infty}{I_N(\a)\over\ln^2N}.$$
\endproclaim

{\bf 5.4.} Observe that Podkorytov in [P2] and Skopina in [Sk1, Sk2] gave
some asymptotic estimates for more general linear means in the cases
which we may treat as "polyhedral" as well. Let
$$\rho(x)=\rho_E (x) = \inf\{\alpha > 0 : {x \over \alpha} \in E \}$$
be the Minkowski functional of a set $E$ and
$$L_N^{\lambda}(f;x)=L_N^{\lambda_E}(f;x)=\sum\limits_{k\in NE}
\lambda\biggl({\rho(k)\over N}\biggr)\hat f(k) e^{ikx}.$$
\proclaim{Theorem 5.4} {\rm ([P2])} Let $E$ be a polyhedron starlike
with respect to the origin, which is an interior point of it, and
$\lambda\in C[0,\infty)$ be supported on $[0,1].$ \medskip\flushpar
{\bf 1.} If the extension of at least one of the faces of the polyhedron
$E$ passes through the origin, then $$\sup\limits_N
||L_N^{\lambda_E}||=\infty$$ and consequently there exists an $f\in
C(\T)$ such that $$\varlimsup\limits_{N\to\infty}|L_N^{\fa}(f;0)|=\infty.$$
{\bf 2.} If the extension of all faces do not pass through the origin,
then the convergence of the integral
$$F_n(\lambda)=\int\limits_{\bold R}|d\hat\lambda(r)|{\ln^{n-1}(2+|r|)
\over 1+|r|}\,dr,$$ where $$d\hat\lambda(r)=\int\limits_0^1 e^{-irt}\,
d\lambda(t),$$ is sufficient for the norms $||L_N^\lambda||$ to be
bounded, and consequently $L_N^\lambda(f;\cdot)$ converge uniformly
to $f$ as $N\to\infty$ for all $f\in C(\T).$   \endproclaim

Some results for "polyhedral" functions
$\lambda$ are obtained in [Sk2, Sk3] in the form like in Theorem 4.1.
In particular, the following  relation holds.
\proclaim{Theorem 5.5} {\rm ([Sk3])} Let $E$ be an $n$-dimensional polyhedron
with rational vertices, starlike with  respect to the origin, and the origin
does not lie at extension of any face of the polyhedron.
Let $\fa (x) = \fa_E.$ Then
$${\Vert  L_{N}^{\lambda}\Vert}_{L_1 (\T) \to L_1 (\T)} = (2\pi)^{-n}
\int\limits_{N\T}
|\faf (x)| \, dx + O( V_{\fa_0} + |\fa_0 (0)|)\ln^{n-1} N . $$ \endproclaim
\bigskip\flushpar
On the base of this theorem, it is possible to find the main term of
$${\Vert L_N^{\fa}\Vert}_{L_1 (\T) \to L_1 (\T)}$$ in the form suitable for
calculations. It is shown in [Sk3], that:
\bigskip
\proclaim\nofrills{ } If $E$ is a convex symmetric $2l$-polygon,
and $\fa_0\in C[0,1]\cap C^1[0,1)$ is such that $\fa_0 (t) \ge 0$,
$\fa_0(1)=0$, $\fa'_0(t)$ and $(t-1)\fa'_0(t)$ is monotone decreasing, then
  $$\Vert L_N^{\fa}\Vert_{L_1 ({\bold T^2})\to L_1 ({\bold T^2})} =
{16\,l\over\pi^4}\, \int\limits_1^N{{\fa_0(1-{1\over x})\,\ln \, x}\over x}\,
dx + O\biggl(\int\limits_1^N {{\fa_0(1-{1\over x})}\over x}\, dx +
\lambda_0(0)\, \biggr) .$$    \endproclaim
\bigskip\flushpar
This allows us to obtain the logarithmic asymptotics providing
some estimates of remainder values. It is found that the constant in the main
term depends on geometric properties of the polyhedron. It is shown in [Sk1]
that the Lebesgue constants grow as $({2\over\pi})^{2n} \ln^n N$ for
parallelepipeds, and as ${2(n+1)\over{\pi^{n+1}}} \ln^n N$ for simplexes.
More precisely, let for $N=0,1,2,...,$ and $0\le
p\le N$ the means $L_N^{\fa}$ are defined as follows:
$$\lambda (x)=\cases 1, & \text{for}\quad x\in(N-p)E,\\
{N+1-\rho(x)\over p+1},&\text{for}\quad x\in NE\setminus(N-p)E,\\
0,&\text{for}\quad x\not\in NE,\endcases$$ where $E$ is the same as in
Theorem 5.5. She rpoved that the norms of such operators are
$$||L_N^{\fa_E}||=(2\pi)^{-n}\int\limits_{\Bbb T^n}|\hat\lambda(x)|\,dx
+\Sigma,$$ where $$|\Sigma|\le C_{P,n}{1\over p+1}(\ln(N+2))^{n-1}.$$

\head 6. "Hyperbolic" results. \endhead

Since the appearance of Babenko's paper \cite{Ba1} interest has continued
in various questions of Approximation Theory and Fourier Analysis in $\R$
connected with studying linear means with harmonics in ``hyperbolic crosses''
$$\Gamma(N,\g)=\{k \in \Zn : h(N,k,\g) = \prod \limits_{j=1}^{n}
(\vert k_j \vert /N)^{\g_j} \le 1, \quad \g_j > 0, j=1,...,n \} .$$
We are interested in the hyperbolic means of Bochner-Riesz type of
order $\a \ge 0$
$$L_{\Gamma(N,\g)}^{\a} : f(x) \mapsto \sum \limits_{k \in \Gamma (N,\g)}
\left(1-h (N,k,\g) \right)_+^{\a} \hat f(k) e^{ik\cdot x} .$$

Hyperbolic Bochner-Riesz means (for the two-dimensional Fourier integrals
with $\g_1 =\g_2 =2$) appeared firstly in the paper of El-Kohen \cite{EK}
in connection with the study of its $L^p$-norms. By the way, his result
was not sharp, and shortly after was strengthened by Carbery \cite{C}.

Hyperbolic partial sums $L_{\Gamma(N,\g)}=L_{\Gamma(N,\g)}^0 $
were investigated separately earlier. The exact degree of growth for them
$\Vert L_{\Gamma(N,\g)} \Vert \asymp N^{\Cp}$ (cf. Theorem 1.1) was
established in the two-dimensional case independently by Belinskii \cite{Be2}
and by A.A. and V.A. Yudins \cite{YY1}, and afterwards was generalized to the
case of arbitrary dimension by Liflyand \cite{L1}.
For $\a>0$ the following relations hold.

\proclaim{Theorem 6.1} {\rm ([L4])}
\medskip
1) For $\a<\Cp \ $ we have  $\Vert L_{\Gamma(N,\g)}^{\a} \Vert
\asymp N^{\Cp-\a}$.
\medskip
2) $\Vert L_{\Gamma(N,\g)}^{\Cp} \Vert = \omega_{n,\g} \ln^n N +
O( \ln^{n-1} N ) .$
\medskip
3) For $\a>\Cp \ $ we have $\Vert L_{\Gamma(N,\g)}^{\a} \Vert =\omega_{n,\g,\a}
\ln^{n-1} N  + O ( \ln^{n-2} N ) .$
\endproclaim\medskip\flushpar
Here and below $\omega$ with subscripts denotes, generally saying, different
constants depending only on the indicated indices.

Observe that the critical order $\Cp$ is the same as in the spherical
case. But if for the values lower than the critical one the orders of growth
of the \Lc\ coincide, the difference between  (3.11) and 2) in Theorem 6.1 is
obvious as well as for orders greater than $\Cp$: in this case the \Lc\
of the usual Bochner-Riesz spherical means are bounded.
The lower estimate in the case 1) follows also from Theorem 3.2.
In order to establish Theorem 6.1, especially 2) and 3), we need the following

\proclaim{Theorem 6.2} {\rm ([L4])} For the norms of operators
$$\bar{L}_{\Gamma(N,\g)}^{\a} :\qquad f(x) \mapsto \sum \limits_{|k_j| \le N,
j=1,\dots,n} \left( 1- h(N,k,\g) \right)_+^{\a} {\hat f}(k) \ek $$
the following asymptotic equality is true
$$\Vert \bar{L}_{\Gamma(N,\g)}^{\a} \Vert = \omega_{n,\g,\a}
\ln^{n-1} N + O ( \ln^{n-2} N ) .$$
\endproclaim
This is a strengthening of Kivinukk's result \cite{Ki}, who was the
first found the effect of influence of the smoothness at the corner points
on the order drop of a logarithmic growth, as compared with the \Lc\
of cubic partial sums, and established double-sided ordinal inequalities.

It should be said that these theorems are proved by step by step passage
from sums to corresponding integrals. This leads to the Fourier transform
of a function generating the method of summability under consideration.

\head 7. Appendix \endhead
In this section we collected results on Lebesgue constants
which are not (explicitly) connected with the Fourier transform
methods. Let us start with one result due to Trigub [T4].

Let $\{e_j\}_{j=1}^n$ be the standard basis in $\R,$ $M_0=
(1,...,n),$ and $q=\sum q_j e_j$ where the $q_j$ are natural numbers
($j\in M_0$); analogously $h=\sum h_j e_j$ where the
$h_j$ are also natural numbers. Set
$$\Delta_{h_j}\lambda_k=\lambda_k-\lambda_{k+h_j e_j}$$
(the difference operator with stepsize $h_j$ in the direction $e_j$) and
$$\Delta_h^q\lambda_k=\biggl(\prod\limits_{j\in M_0}\Delta^{q_j}_{h_j}
\biggr)\lambda_k$$ ("mixed" difference in the direction of all axes).
\proclaim{Theorem 7.1} For every $p\in[1,2)$ and $q,$ there exists a constant
$C,$ depending only on $p,$ $q$ and $m,$ such that
$$\int\limits_{\T}|\sum\limits_{-N_j\le k_j\le N_j}\lambda_k e^{ikx}|^p
\,dx\le C\prod\limits_j(N_j+1)^{{p-2\over 2}}\sum\limits_{0\le s_j\le
[\log_2(N_j+1)]} 2^{{2-p\over 2}\sum\limits_j s_j}(\sum\limits_k|\Delta_h^q
\lambda_k|^2)^{{p\over 2}}$$
where $\lambda_k$ is taken to equal $0$ for $k_j\ne[-N_j,N_j]$ in the sum
$\sum\limits_k,$ while $h=h(s,N)$ is defined by the following conditions
$${N_j+1\over 3\cdot 2^{s_j}}\le h_j\le {5(N_j+1)\over 6\cdot 2^{s_j}},
\qquad\qquad {N_j+1\over 3\cdot 2^{s_j}}\le h_j\le {N_j+1\over 2^{s_j}}$$
according as $s_j<[\log_2(N_j+1)]$ or $s_j=[\log_2(N_j+1)].$\endproclaim
In several corollaries sufficient conditions are given to ensure that the
Lebesgue constants have a given rate of growth. This is done in terms of
smoothness of a function generating the sequence $\{\lambda_k\},$
namely $\lambda_k=\lambda_{N,k}=\lambda({k_1\over N_1},...,{k_n\over N_n}).$

Let us go on to results in which the operator of taking of
partial sums is unbounded. The first one is due to Belinskii [Be5]
(see also [MP]).

Let $l_1, l_2,...,l_k,$ with $k<n,$ be a family of linear independent
vectors in $\R.$ Consider the sets
$$P_j=\{m\in\bold Z^n: |l_j m|\le 1\}$$ and put $$P=\bigcap_{j=1}^k P_j$$
and $$P_0=\bigcap_{j=1}^k\{x\in\R: l_j x=0\}.$$
\proclaim{Theorem 7.2} {\bf 1.} If in $P_0$ there exists a sublattice
of $\bold Z^n$ of dimension $n-k,$ then
$$||S_{NP}||\asymp \ln^k N\qquad\text{as}\quad N\to\infty.$$
{\bf 2.} If there not exist anyone such sublattice in $P_0$ then the
operator $S_{NP}$ is unbounded for each $N>0.$ \endproclaim
\remark{Remark} If in the first case of Theorem 7.2 for a set considered
an asymptotics is proved (see Section 5 or also [Ku2]) then one gets
the asymptotics in Theorem 7.2 as well.  \endremark

The next theorem, due to Belinskii and Liflyand [BL3],
is of the same nature but deals with hyperbolic partial sums
(see also Section 6). Let $L_j(x)=l_{j1}x_1+...+l_{jn}x_n,$
$j=1,2,...,n,$ be linear forms with nonsingular coefficient matrix
$\Lambda=\{l_{jk}\},$ $1\le j,k\le n,$ $\det \Lambda\ne 0,$ and
$$H=\{x\in\R: \prod\limits_{j=1}^n |L_j(x)|\le 1\}.$$
We call the matrix $\Lambda$ {\it rational} if each row of it consists
of integers, up to a common factor. In the contrary case, the matrix
is said to be {\it irrational.}
\proclaim{Theorem 7.3} {\bf 1.} If the matrix $\Lambda$ is {\it rational,} then
$$||S_{NH}||\asymp N^{{n-1\over 2}}.$$ {\bf 2.} If $\Lambda$ is
{\it irrational,} then there exists an integer $N_0$ such that the operator
$S_{NH}$ is unbounded for all $N>N_0.$ \endproclaim
In both theorems the second parts are proved by using the following
result from [Ru1] (see also [Ru2], Th.3.1.3):

\flushpar  {\sl If the operator
of taking partial sums with respect to some dilation of a given set
is bounded than this set may be represented as a finite union of cosets
of discrete subgroups of the lattice} $\bold Z^n.$

\flushpar  To get a contradiction with this statement some theorems
from Geometric Number Theory are used (see [Ca]).

For $n=2$ Theorem 7.3 was earlier obtained by Belinskii (see [Be3, Be4]).
The proof was essentially two-dimensional and relied on some other results
in Number Theory, in particular an exact value $N_0$ was indicated.
Nevertheless even in this case the second part of Theorem 7.3 cannot be
established for more general sets $H=\{\prod |L_j|^{\gamma_j}\le 1\},$
since there are no corresponding results in Number Theory.

The following result obtained by A. Yudin and V. Yudin [YY1] is closely
connected with the result of Podkorytov given above in Theorem 1.2.
The latter is even easier for proving the estimate from above for the Lebesgue
constants of hyperbolic partial sums (see Section 6) though it was not
applied to this and the author knows the proof from Podkorytov's private
communication. In contrary, in [YY1] what follows was
invented especially for "hyperbolic" estimates. Thus, let $U\subset\bold Z^n$
be a bounded set and $t\in\bold Z^n.$ We set
$$U_t=\{k\in\bold Z^n: k-t\in U\}$$ and
$$\omega(t,U)=2|U|-|U\cap U_t|-|U\cap U_{-t}|,$$
where $|U|$ denotes the number of points in $U.$
\proclaim{Theorem 7.4} Let numbers $L_1\le L^2$ be such that
$$\omega(he_j,U)\le L_1 h\qquad\text{and}\qquad\omega(he_r,U)\le L_2 h$$
for some natural numbers $r,j\in M_0,$ where $r\ne j,$ and every natural
number $h.$ Then $$||S_U||\le{1\over 2}({L_1\over 2})^{{1\over 2}}
\log_2{L_2\over L_1}+{3\over 2-\sqrt{2}}L_1^{{1\over 2}}.$$\endproclaim

The Lebesgue constants of step-wise hyperbolic crosses have been considered
in many papers, together with various applications of such estimates.
These problems were discussed by Temlyakov [Tm], Galeev [Ga1, Ga2],
and Belinskii. For example it was shown by Belinskii in [Be6] that
if $H_N$ is defined as
$$H_N=\bigcup\{m\in\bold Z^n: 2^{s_j}\le|m_j|<2^{s_j+1}\}$$
for $s\in\bold Z^n\cap[0,\infty)^n$ such that $0\le s_1+...+s_n\le N,$ and
$N=1,2,...,$ then $$||S_{H_N}||\asymp N^{n+{n-1\over 2}}$$ as $N\to\infty.$

We now need to introduce some new notation to be able to formulate
further results due to M. Dyachenko. Let $A_2$ be the class of bounded
sets $U\subset\bold Z^n$ such that if $m\in U,$ then
$$\bold Z^n\bigcap\prod\limits_{j=1}^n [\min(m_j,0),\max(m_j,0)]\subset U,$$
and let us define $A_1$ by
$$A_1=\{U\bigcap(0,\infty)^n,\quad\text{where}\quad U\in A_2\}.$$
We also define $M_1$ as the class of $n$-dimensional sequences
$$a=\{a_m\}=\{a_{m_1,...,m_n}\}_{m_1,...,m_n=1}^\infty$$
such that $1\le k_j\le m_j$ implies that $a_k\ge a_m\ge 0.$ Set
$$\prod(x)=\prod\limits_{j=1}^n(|x_j|+1),$$ and it is then possible
to give the following assertions:
\proclaim{Theorem 7.5} Given $U\in A_1$ or $U\in A_2$ and a number
$p\in[1,{2n\over n+1}),$ then
$$||S_U||_{L^p}\le C_{p,n}\max\limits_{m\in U}
\biggl(\prod(m)\biggr)^{{n-1\over 2n}}.$$
\endproclaim\flushpar
We note that Theorem 7.5 yields the upper bound of Theorem 6.1 in
the case 1) for $\alpha=0$ and $\gamma_1=...=\gamma_n=1.$
\proclaim{Corollary 7.1} The following inequalities are satisfied under
the hypotheses of Theorem 7.5:
$$||S_U||_{L^p}\le C_{p,m}|U|^{{n-1\over 2n}}$$  and
$$||S_U||_{L^p}\le C_{p,m}\biggl(\sum\limits_{m\in U}(\prod(m))^
{{2n\over n+1}-2}\biggr)^{{n+1\over 2n}}.$$  \endproclaim
The first inequality of Corollary 7.1 was proved in [Dy2] for $p=1$
and $U\in A_1$ with the constant $C_{1,n}=50n^3.$ In [Dy1] such an
estimate was obtained with an additional logarithmic factor. Some
other estimates for $p>1$ as well as some open problems can be found
in the survey [Dy3], Sect.3.

Mention also results of Ustina on the two-dimensional Hausdorff method
(see [U]).

\enddocument